\documentclass[letterpaper,twocolumn,10pt]{IEEEtran}
\usepackage{multicol}
\usepackage{multirow}
\usepackage{lipsum}
\usepackage{booktabs}
\usepackage{amsfonts}
\usepackage{amsmath,cases}
\usepackage{amssymb}
\usepackage{graphicx}
\usepackage{cite}
\usepackage{epsfig}
\usepackage[T1]{fontenc}
\usepackage[lowtilde]{url}
\usepackage{color}
\usepackage{algorithm}
\usepackage{algorithmic}
\usepackage{epstopdf}
\usepackage{color}
\usepackage [english]{babel}
\usepackage [autostyle, english = american]{csquotes}
\MakeOuterQuote{"}
\usepackage{color}
\usepackage{amssymb}
\usepackage{bm}
\usepackage{caption}
\usepackage{subcaption}

\usepackage{float}
\usepackage{multirow}
\usepackage{hhline}

\usepackage{algorithm,algorithmic}
\usepackage{array}
\ifCLASSOPTIONcompsoc
\usepackage[caption=false,font=footnotesize,labelfont=sf,textfont=sf]{subfig}
\else
\usepackage[caption=false,font=footnotesize]{subfig}
\fi
\usepackage{fixltx2e}
\usepackage{amsthm}
\theoremstyle{plain}

\newcommand{\Tr}{{\sf Trace}}

\newcommand{\al}{\alpha}

\newcommand{\mc}{\mathcal}
\newcommand{\mb}{\mathbb}

\newcommand{\tb}{\textbf}
\newcommand{\nn}{\nonumber}
\newcommand{\bea}{\begin{eqnarray}}
\newcommand{\eea}{\end{eqnarray}}
\newcommand{\beq}{\begin{equation}}
\newcommand{\eeq}{\end{equation}}

\newcommand{\bet}{\beta}

\newtheorem{ex}{Example}\newcommand{\Ex}{\begin{ex}\rm}
	\newcommand{\eex}{\end{ex}}

\begin{document}

\title{Efficient tensor completion for color image and video recovery: Low-rank tensor train}
\author{Johann A. Bengua$^1$, Ho N. Phien$^1$, Hoang D. Tuan$^1$\thanks{$^1$Faculty of Engineering and Information Technology,	University of Technology Sydney, Ultimo, NSW 2007, Australia; Email: johann.a.bengua@student.uts.edu.au, ngocphien.ho@uts.edu.au, tuan.hoang@uts.edu.au} and Minh N. Do$^2$\thanks{
$^2$Department of Electrical and Computer Engineering and the Coordinated Science Laboratory, University of Illinois at Urbana-Champaign, Urbana, IL 61801 USA; Email: minhdo@illinois.edu }}

\maketitle

\vspace*{-1.0cm}

\begin{abstract}
This paper proposes {\color{black} a novel approach to tensor completion, which
recovers missing entries of data represented by tensors.} The approach is based on the tensor train (TT) rank,
{\color{black} which  is able to capture hidden information from tensors thanks to its definition from a well-balanced matricization scheme.} {\color{black}Accordingly, new optimization formulations for tensor completion are proposed as well as two new algorithms for their solution}. The first one called simple low-rank tensor completion via tensor train (SiLRTC-TT) is intimately related to minimizing a nuclear norm based on TT rank. The second one is from a multilinear matrix factorization model to approximate the TT rank of a tensor, and is called tensor completion by parallel matrix factorization via tensor train (TMac-TT). {\color{black}A tensor augmentation scheme of transforming a low-order tensor to higher-orders is also proposed to enhance the effectiveness of SiLRTC-TT and TMac-TT.}
{\color{black}Simulation results for color image and video recovery show the clear advantage of our method over all other methods.}
\end{abstract}

	\title{}	
	\author{
		\IEEEcompsocitemizethanks{\IEEEcompsocthanksitem The authors are from the Faculty of Engineering and Information Technology, University of Technology Sydney, Ultimo, NSW 2007, Australia.\protect\\ E-mail: ngocphien.ho@uts.edu.au, tuan.hoang@uts.edu.au.\protect}\\}	
\section{Introduction}
Tensors are multi-dimensional arrays, which are higher-order generalizations of matrices and vectors \cite{Tamara_2009}. Tensors provide a natural way to represent multidimensional data whose entries are indexed by several continuous or discrete variables.  Employing tensors and their decompositions to process data has become increasingly popular since
\cite{Vasilescu_2003,Sun_2005,Franz_2009}. For instance, a color image is a third-order tensor defined by two indices for spatial variables and one index for color mode. A video comprised of color images is a fourth-order tensor with an additional index for a temporal variable. Residing  in extremely high-dimensional data spaces, the tensors in practical applications are nevertheless often of
\emph{low-rank} \cite{Tamara_2009}. Consequently, they can be effectively  projected to much smaller subspaces through underlying decompositions such as the CANDECOMP/PARAFAC (CP)\cite{Carroll_1970,Harshman_1970}, Tucker \cite{Tucker_1966} and tensor train (TT) \cite{Oseledets_2011} or matrix product state (MPS)\cite{Fannes_1992,Klumper_1993,PerezGarcia_2007}.

Motivated by the success of low rank matrix completion (LRMC) \cite{Cai_2010,Ma_2009,Recht_2010}, recent effort has been made to extend its concept to low rank tensor completion (LRTC). In fact, LRTC {\color{black}has found applications} in computer vision and graphics, signal processing and machine learning \cite{Ji2013,Signoretto_2010,Signoretto_2011,Gandy_2011,Tomioka_2011,Tan_2014,Xu}. The common target is to recover  missing entries of a tensor from its partially observed entities \cite{Bertalmio_2000,Komodakis_2006,Korah_2007}. LRTC remains a grand challenge due to the fact that {\color{black} computation for the tensor rank, defined as CP rank, is already an NP-hard problem
\cite{Tamara_2009}.}
 There have been {\color{black} attempts in approaching} LRTC via \emph{Tucker rank} \cite{Ji2013,Gandy_2011,Xu}. A {\color{black}conceptual} drawback of Tucker rank is that its components are ranks of matrices constructed based on an unbalanced matricization scheme (one mode versus the rest). The upper bound of each individual rank is often small and may not be suitable for describing global information of the tensor. In addition, the matrix rank minimizations is only efficient when the matrix is more balanced. As the rank of a matrix is not more than $\min\{n,m\}$, where $m$ and $n$ are the number of rows and columns of the matrix, respectively, the high ratio $\max\{m,n\}/\min\{m,n\}$ would effectively rule out the need of matrix rank minimization.
{\color{black} It is not surprising for present state-of-the-art LRMC methods
\cite{Bach_2008,Candes_2009, Cai_2010,Ma_2009} to implicitly assume that the considered matrices are balanced.}

{\color{black}Another type of tensor rank is the \emph{TT rank}, which constitutes of ranks of matrices formed by a well-balanced matricization scheme, i.e. matricize the tensor along permutations of modes. TT rank was defined in \cite{Oseledets_2011}, yet low rank tensor analysis via TT rank can be seen in earlier work in physics, specifically in simulations of quantum dynamics \cite{PhysRevLett.91.147902, Vidal_2004}.  Realizing the computational efficiency of low TT rank tensors,
there has been numerous works in applying it to numerical linear algebra \cite{otz-teninv-2011, corona2015, doi:10.1137/110833142}.
Low TT rank tensors were used for the singular value decomposition (SVD) of large-scale matrices in \cite{Mach2013,doi:10.1137/140983410}. The alternating least squares (ALS) algorithms for  tensor approximation \cite{Holtz_2012,Grasedyck_2015} are also used for solutions of linear equations and eigenvector/eigenvalue approximation.
In \cite{Da_Silva_2015,Rauhut_2015}, low TT rank tensors were also used in implementing
the steepest descent iteration for large scale least squares problems. The common assumption in all these works is that
all the used tensors during the computation processes are of low TT rank for computational practicability. How
low TT rank tensors are relevant to real-world problems was not really their concern. Applications of the TT decomposition to fields outside of mathematics and physics has rarely been seen, with only a recent application of TT to machine learning \cite{7207289}.
As mentioned above, color image and video are perfect examples of tensors, so their completion can be formulated
as tensor completion problems. However, it is still not known if TT rank-based completion is useful for practical solutions.
The main purpose of this paper is to show that TT rank is the right approach for LRTC, which can be addressed
by TT rank-based optimization. The paper contribution is as follows:
\begin{enumerate}
\item Using the concept of von Neumann entropy in quantum information theory \cite{Nielsen_2009},
we show that Tucker rank does not capture the global correlation of the tensor entries and thus is hardly ideal for
LRTC. Since TT rank constitutes of ranks of matrices formed by a well-balanced matricization scheme,
it is capable of capturing the global correlation of the tensor entries and is thus a promising tool for LRTC.
\item We show that unlike Tucker rank, which is often low and not interesting for optimization,
TT rank optimization is a tractable formulation for LRTC.
Two new algorithms are introduced to address the TT rank optimization based LRTC problems.
The first algorithm called simple low-rank tensor completion via tensor train (SiLRTC-TT)
  solves an optimization problem based on the \emph{TT nuclear norm}. The second algorithm called tensor completion by parallel matrix factorization via tensor train (TMac-TT) uses a mutilinear matrix factorization model to approximate the TT rank of a tensor, bypassing the computationally expensive SVD. Avoiding the direct TT decomposition enables the proposed algorithms to outperform other start-of-the-art tensor completion algorithms.
\item We also introduce a novel technique called \emph{ket augmentation} (KA) to represent a low-order tensor by a higher-order tensor without changing the total number of entries. The KA scheme  provides a perfect means to obtain a higher-order tensor representation of visual data by maximally exploring the potential of TT rank-based optimization for color image and video completion. TMac-TT especially performs well in recovering videos with $95\%$ missing entries.
\end{enumerate}}
The rest of the paper is organized as follows. Section \ref{sec1} introduces notation and a brief review of tensor  decompositions. In Section \ref{sec2}, the conventional formulation of LRTC is reviewed, the advantage of TT rank over Tucker rank in terms of global correlations is discussed, then the proposed reformulations of LRTC in the concept of TT rank.  Section \ref{propalg} introduces two algorithms to solve the LRTC problems based on TT rank, followed by a discussion of computational complexity. Subsequently, the tensor augmentation scheme known as KA is proposed in Section \ref{sec3}. Section \ref{sec4} provides experimental results and finally, we conclude our work in Section \ref{sec5}.

\section{Tensor ranks \label{sec1}}
Some mathematical notations and preliminaries of tensors are adopted from \cite{Tamara_2009}. A \textit{tensor} is a multi-dimensional array and its \textit{order} or \emph{mode} is the number of its dimensions. Scalars are zero-order tensors denoted by lowercase letters $(x, y, z,\ldots)$. Vectors and matrices are the first- and second-order tensors which are denoted by boldface lowercase letters (\textbf{x}, \textbf{y}, \textbf{z},\ldots) and capital letters $(X, Y, Z,\ldots)$, respectively. A higher-order tensor (the tensor of order three or above) are denoted by calligraphic letters $(\mc{X}, \mc{Y}, \mc{Z},\ldots)$.

An \textit{N}th-order tensor is denoted as $\mc{X}\in\mathbb{R}^{I_1\times I_2\times\cdots\times I_N}$ where $I_k$, $k=1,\ldots, N$ is the dimension corresponding to mode $k$. The elements of $\mc{X}$ are denoted as $x_{i_1\cdots i_k\cdots i_N}$, where $1\leq i_k\leq I_k$, $k=1,\ldots, N$. The Frobenius norm of $\mc{X}$ is $||\mc{X}||_F = \sqrt{\sum_{i_1}\sum_{i_2}\cdots\sum_{i_N}x^2_{i_1i_2\cdots i_{N}}}$.

A mode-$n$ fiber of a tensor $\mc{X}\in\mathbb{R}^{I_1\times I_2\times\cdots\times I_N}$ is a vector defined by fixing all indices but $i_n$ and denoted by \text{\bf x}$_{i_1\ldots i_{n-1}:i_{n+1}\ldots i_{N}}$.

Mode-$n$ matricization (also known as mode-$n$ unfolding or flattening) of a tensor $\mc{X}\in\mathbb{R}^{I_1\times I_2\times\cdots\times I_N}$ is the process of unfolding or reshaping the tensor into a matrix $X_{(n)}\in\mathbb{R}^{I_n\times (I_1\cdots I_{k-1}I_{k+1}\cdots I_N)}$ by rearranging the mode-$n$ fibers to be the columns of the resulting matrix.  Tensor element $(i_1,\ldots, i_{n-1},i_n,i_{n+1},\ldots, i_{N})$ maps to matrix element $(i_n,j)$ such that
\bea
j=1+\sum_{k=1,k\neq n}^{N}(i_k-1)J_k~~\text{with}~~J_k=\prod_{m=1, m\neq n}^{k-1}I_m.
\label{indexj}
\eea
The vector $\tb{r} = (r_1, r_2,\ldots, r_N)$, where $r_n$ is the rank of the corresponding matrix $X_{(n)}$ denoted as $r_n = \text{rank}(X_{(n)})$, is called as the \emph{Tucker rank} of the tensor $\mc{X}$. {\color{black}It is obvious that $\text{rank}(X_{(n)})\leq I_n$.}

Using Vidal's decomposition \cite{Vidal_2004}, $\mc{X}$ can be {\color{black}represented by a sequence of connected low-order
tensors   in the form}
\bea
\mc{X}
&=&\sum_{i_1,\ldots,i_N}\Gamma^{[1]}_{i_{1}}\lambda^{[1]}\cdots\lambda^{[N-1]}\Gamma^{[N]}_{i_{N}}\tb{e}_{i_1}\otimes\cdots\otimes\tb{e}_{i_N},~~~
\label{tt2}
\eea
where for $k=1,\ldots,N$, $\Gamma^{[k]}_{i_k}$ is an $r_{k-1}\times r_k$ matrix and $\lambda^{[k]}$ is the $r_{k}\times r_{k}$ diagonal singular matrix, $r_0=r_{N+1}=1$. For every $k$, the following orthogonal conditions are fulfilled:
\bea
\sum_{i_k=1}^{I_k}\Gamma^{[k]}_{i_{k}}\lambda^{[k]}(\Gamma^{[k]}_{i_{k}}\lambda^{[k]})^{T} &=& \mb{I}^{[k-1]},\\
\sum_{i_k=1}^{I_k}(\lambda^{[k-1]}\Gamma^{[k]}_{i_{k}})^{T}\lambda^{[k-1]}\Gamma^{[k]}_{i_{k}} &=& \mb{I}^{[k]},
\eea
where $\mb{I}^{[k-1]}$ and $\mb{I}^{[k]}$ are the identity matrices of sizes $r_{k-1}\times r_{k-1}$ and $r_{k}\times r_{k}$, respectively. The \emph{TT rank} of the tensor is simply defined as $\tb{r} = (r_1,r_2,\ldots, r_{N-1})$, and can be determined directly via the singular matrices $\lambda^{[k]}$. Specifically, to  determine $r_k$, rewrite (\ref{tt2}) as
\bea
\mc{X}
&=&\sum_{i_1,i_2\ldots,i_N}\tb{u}^{[1\cdots k]i_1\cdots i_k}\lambda^{[k]}\tb{v}^{[k+1\cdots N]i_{k+1}\cdots i_N},
\label{tt3}
\eea
where
\bea
\tb{u}^{[1\cdots k]i_1\cdots i_k}&=&\Gamma^{[1]}_{i_{1}}\lambda^{[1]}\cdots\Gamma^{[k]}_{i_{k}}\otimes_{l=1}^{k}\tb{e}_{i_l},
\eea
and
\bea
\tb{v}^{[k+1\cdots N]i_{k+1}\cdots i_N}&=&\Gamma^{[k+1]}_{i_{k+1}}\lambda^{[k+1]}\cdots\Gamma^{[N]}_{i_{N}}\otimes_{l=k+1}^{N}\tb{e}_{i_l}.~~~~~~
\eea
We can also rewrite (\ref{tt3}) in terms of the matrix form of an SVD as
\bea
X_{[k]}&=& U\lambda^{[k]}V^{T},
\label{tt4}
\eea
where $X_{[k]}\in\mb{R}^{m\times n}$ ($m = \prod_{l=1}^{k}I_l, n=\prod_{l=k+1}^{N}I_l$) is  the \emph{mode-$(1,2,\ldots, k)$ matricization} of the tensor $\mc{X}$ \cite{Oseledets_2011}, $U\in\mb{R}^{m\times r_k}$ and $V\in\mb{R}^{n\times r_k}$ are orthogonal matrices. The rank of $X_{[k]}$ is $r_k$, which is defined as the number of nonvanishing singular values of $\lambda^{[k]}$.

Since matrix $X_{[k]}$ is obtained by matricizing along {$k$ modes,
its rank $r_k$  is bounded by $\min(\prod_{l=1}^{k}I_l,\prod_{l=k+1}^{N}I_l)$.

\section{Tensor completion \label{sec2}}
{\color{black}This section firstly revisits the conventional formulation of LRTC based on the Tucker rank. Then,
we propose a new approach to LRTC via TT rank optimization, which leads to two new optimization formulations, one based on nuclear norm minimization, and the other on multilinear matrix factorization.}
\subsection{Conventional tensor completion}
{\color{black}As tensor completion is fundamentally based on matrix completion, we give
an overview of the latter prior its introduction.} {\color{black}Recovering} missing entries of a  matrix  $T\in\mb{R}^{m\times n}$ from its partially known entries given by a subset $\Omega$ can be studied via the well-known matrix-rank optimization problem \cite{fazel2002matrix,Kurucz_2007}:
\bea
\begin{aligned}
	& \underset{X}{\text{min}}&&\text{rank}(X)~~~~ \text{s.t.}&& X_{\Omega} = T_{\Omega}.
\end{aligned}
\label{eq1}
\eea
The missing entries of $X$ are completed such that the rank of $X$ is as small as possible, i.e
the vector $(\lambda_1,....,\lambda_{\min\{m,n\}})$ of the singular values $\lambda_k$ of $X$ is as sparse as possible. The sparsity
of $(\lambda_1,....,\lambda_{\min\{m,n\}})$ leads to the effective representation of $X$ for accurate completion.
Due to the combinational nature of the function $\text{rank}(\cdot)$, the problem (\ref{eq1}), however, is NP-hard.
For the nuclear norm $||X||_{*}=\sum_{k=1}^{\min\{m,n\}}\lambda_k$, the following convex
$\ell^1$ optimization problem in $(\lambda_1,....,\lambda_{\min\{m,n\}})$ has been proved the most effective surrogate for (\ref{eq1})  \cite{Bach_2008,Cai_2010,Ma_2009}:
\bea
\begin{aligned}
	& \underset{X}{\text{min}}&&||X||_{*}~~~~ \text{s.t.}&& X_{\Omega} = T_{\Omega}.
\end{aligned}
\label{eq2}
\eea
It should be emphasized that the formulation (\ref{eq1}) is efficient only when $X$ is
balanced (square), i.e. $m\approx n$.  It is likely that $\mbox{rank}(X)\approx m$ for unbalanced $X$ with $m\ll n$, i.e.
there is not much difference between the optimal value of (\ref{eq1}) and its upper bound $m$, under which
{\color{black}rank optimization problem} (\ref{eq1}) is not interesting.
More importantly, one needs at least $Cn^{6/5}\mbox{rank}(X)\log n\approx Cn^{6/5}m\log n$ sampled entries \cite{Candes_2009} with a positive constant $C$ to successfully complete $X$, which is almost the total $nm$ entries of $X$.

{\color{black}Completing} an $N$th-order tensor $\mc{T}\in\mathbb{R}^{I_1\times I_2\cdots\times I_{N}}$ from its known entries given by an index set $\Omega$ is {\color{black}formulated by} the following
{\color{black}Tucker rank} optimization problem \cite{Ji2013,Gandy_2011,Tan_2014,Xu}:
\bea
\begin{aligned}
	& \underset{X_{(k)}}{\text{min}}&&\sum_{k=1}^{N}\al_{k}\text{rank}(X_{(k)})~~~~ \text{s.t.}&& \mc{X}_{\Omega} = \mc{T}_{\Omega}.
\end{aligned}
\label{eq4}
\eea
where $\{\al_k\}_{k=1}^{N}$ are defined as weights fulfilling the condition $\sum_{k=1}^{N}\al_k=1$,
{\color{black} which is then addressed  by}  the following $\ell^{1}$ optimization problem \cite{Ji2013}:
\bea
\begin{aligned}
	& \underset{X_{(k)}}{\text{min}}&&\sum_{k=1}^{N}\al_{k}||X_{(k)}||_{*}~~~~ \text{s.t.}&& \mc{X}_{\Omega} = \mc{T}_{\Omega}.
\end{aligned}
\label{eq4_1}
\eea
{\color{black} Each matrix $X_{(k)}$ in (\ref{eq4}) is obtained by matricizing the tensor along one single mode and
 thus is highly unbalanced. For instance, when all the modes have the same dimension ($I_1=\cdots=I_N\equiv I$),
its dimension is $I\times I^{N-1}$. {\color{black}As a consequence, its rank is low, which makes the matrix rank optimization}
formulation (\ref{eq4}) less efficient for completing $\mc{T}$. Moreover, as analyzed above, it also makes
the $\ell^{1}$ optimization problem (\ref{eq4_1}) not efficient in addressing the rank optimization problem
(\ref{eq4}).}

{\color{black}In the remainder of this subsection we show that $\mbox{rank}(X_{(k)})$ is not an appropriate means for capturing the \emph{global correlation} of a tensor as
it provides only the mean of the correlation between a single mode (rather than a few modes) and the rest of the tensor.

Firstly, normalize $\mc{X}$ ($||\mc{X}||_F = 1$) and represent it as:
\bea
\mc{X} = \sum_{i_1,i_2\ldots,i_N}x_{i_1i_2\cdots i_N}\tb{e}_{i_1}\otimes\tb{e}_{i_2}\cdots\otimes\tb{e}_{i_N},
\label{eq_state}
\eea
where "$\otimes$" denotes a tensor product \cite{Tamara_2009},  $\tb{e}_{i_k}\in\mathbb{R}^{I_k}$ form an orthonormal
basis in $\mathbb{R}^{I_k}$ for each $k=1,\ldots,N$. Applying mode-$k$ matricization of $\mc{X}$ results in $X_{(k)}$ representing a pure state of a composite system $AB$ in the space $\mc{H}_{AB}\in\mb{R}^{m\times n}$,
which is a tensor product of two subspaces $\mc{H}_{A}\in\mb{R}^{m}$ and $\mc{H}_{B}\in\mb{R}^{n}$ of dimensions $m = I_k$ and $n=\prod\limits_{l=1,l\neq k}^{N}I_l$, respectively. The subsystems $A$ and $B$ are
seen as two \emph{contigous partitions} consisting of mode $k$ and all other  modes of the tensor, respectively.
It follows from (\ref{eq_state}) that
\bea
X_{(k)} = \sum_{i_k,j}x_{i_kj}\tb{e}_{i_k}\otimes\tb{e}_{j},
\label{eq_state1}
\eea
where the new index $j$ is defined as in (\ref{indexj}), $\tb{e}_{j} = \otimes_{l=1,l\neq k}^{N}\tb{e}_{i_l}\in\mb{R}^{n}$. According to the Schmidt decomposition \cite{Nielsen_2009}, there exist orthonormal bases $\{\tb{u}^{A}_{l}\}$ in $\mc{H}_{A}$ and $\{\tb{v}^{B}_{l}\}$ in $\mc{H}_{B}$ such that,
\bea
X_{(k)} = \sum_{l=1}^{r_k}\lambda_l\tb{u}^{A}_{l}\otimes\tb{v}^{B}_{l},
\eea
where $r_k$ is the rank of $X_{(k)}$, $\lambda_l$ are nonvanishing singular values satisfying $\sum_{l=1}^{r_k}\lambda^{2}_l=1$, $\{\tb{u}^{A}_{l}\}$ and $\{\tb{v}^{B}_{l}\}$ are orthonormal bases. The correlation between two subsystems $A$ and $B$ can be studied via von Neumann entropy defined as \cite{Nielsen_2009}:
\bea
S^A=-\Tr(\rho^A\log_2(\rho^A)),
\label{en1}
\eea
where $\rho^A$ is called the \emph{reduced density matrix operator} of the composite system and computed by taking the partial trace of the density matrix $\rho^{AB}$ with respect to $B$. Specifically, we have
\bea
\rho^{AB} &=& X_{(k)}\otimes(X_{(k)})^T\nn\\
&=&\Big(\sum_{l=1}^{r_k}\lambda_l\tb{u}^{A}_{l}\otimes\tb{v}^{B}_{l}\Big)\otimes\Big(\sum_{j=1}^{r_k}\lambda_j\tb{u}^{A}_{j}\otimes\tb{v}^{B}_{j}\Big)^T.
\eea
Then $\rho^A$ is computed as
\bea
\rho^A &=& \Tr_B(\rho^{AB})\nn\\
&=&\sum_{l=1}^{r_k}\lambda^{2}_l\tb{u}^{A}_{l}\otimes(\tb{u}^{A}_{l})^{T},
\label{en2}
\eea
Substituting (\ref{en2}) to (\ref{en1}) yields
\bea
S^A= -\sum_{l=1}^{r_k}\lambda^{2}_l\log_{2}\lambda^{2}_l.
\eea
Similarly,
\bea
S^B&=&-\Tr(\rho^B\log_2(\rho^B))\nn\\
&=&-\sum_{l=1}^{r_k}\lambda^{2}_l\log_{2}\lambda^{2}_l,
\eea
which is the same with $S^A$, simply $S^A=S^B=S$. This entropy reflects the correlation or \emph{degree of entanglement} between  subsystem $A$ and its complement $B$ \cite{Bennett_1996}. It is bounded by $0\leq S\leq\log_{2} r_k$. Obviously, there is no correlation between subsystems
$A$ and $B$ whenever $S=0$ (where $\lambda_1 = 1$ and the other singular values are zeros).  There exists correlation between
 subsystems $A$ and $B$ whenever $S\neq 0$  with its maxima  $S=\log_{2} r_k$ (when $\lambda_1 = \cdots = \lambda_{r_k} = 1/\sqrt{r_k}$).
 Furthermore, if the singular values decay significantly, e.g. exponential decay, we can also keep a few $r_k$
largest singular values of $\lambda$ without considerably losing accuracy in quantifying the amount of correlation between the subsystems. Then $r_k$ is referred to as the approximate \emph{low rank} of the matrix $X_{(k)}$, which also
 means that the amount of correlation between two subsystems $A$ (of mode $k$) and $B$ (of other modes) is small.
 On the contrary, if two subsystems $A$ and $B$ are highly correlated, i.e. the singular values decay very slowly,
 then $r_k$ needs to be as large as possible.

From the above analysis, we see that the rank $r_k$ of $X_{(k)}$ is only capable of capturing the correlation between one mode $k$ and the others. Hence, the problem (\ref{eq4}) does not take into account the correlation between a few modes and the rest of the tensor, and thus may not be sufficient for completing high order tensors ($N>3$). To overcome this weakness, in the next subsection, we will approach LRTC problems optimizing \emph{TT rank}, which is defined by more balanced matrices and is able to capture the hidden correlation between the modes of the tensor more effectively.}
\subsection{\color{black}Tensor completion by TT rank optimization}}
A new {\color{black}approach to  the LRTC problem in (\ref{eq4}) is to address it by the following TT rank optimization}
\bea
\begin{aligned}
	& \underset{X_{[k]}}{\text{min}}&&\sum_{k=1}^{N-1}\al_{k}\text{rank}(X_{[k]})~~~~ \text{s.t.}&& \mc{X}_{\Omega} = \mc{T}_{\Omega},
\end{aligned}
\label{eq5}
\eea
where $\al_{k}$ denotes the weight that the TT rank of the matrix $X_{[k]}$ contributes to, with the condition $\sum_{k=1}^{N-1}\al_k=1$. {\color{black}Recall that $X_{[k]}$ is obtained by matricizing along $k$ modes and
thus its rank {\color{black}captures the correlation between $k$ modes and the other $N-k$ modes}. Therefore,
$(\mbox{rank}(X_{[1]}), \mbox{rank}(X_{[2]}),...,\mbox{rank}(X_{[N]}))$ provides a much better means to capture
the global information of the tensor.}

As the problem (\ref{eq5}) is still difficult to handle as $\text{rank}(\cdot)$ is presumably hard. Therefore, from (\ref{eq5}), we propose the following two problems.

The first one based on the so-called \emph{TT nuclear norm}, defined as
\bea
||\mc{X}||_{*}&=&\sum_{k=1}^{N-1}\al_{k}||X_{[k]}||_{*},
\label{TTnorm}
\eea
is given by
\bea
\begin{aligned}
	& \underset{\mc{X}}{\text{min}}&&\sum_{k=1}^{N-1}\al_{k}||X_{[k]}||_{*}~~~~ \text{s.t.}&& \mc{X}_{\Omega} = \mc{T}_{\Omega},
\end{aligned}
\label{eq6}
\eea
{\color{black}The concerned matrices in (\ref{eq6}) are much more balanced than their counterparts in (\ref{eq4_1}). As a result, the $\ell^1$ optimization problem (\ref{eq6}) provides an effective means for
the matrix rank optimization problem (\ref{eq5}). }\\
{\color{black}A particular case of  (\ref{eq6}) is the square model \cite{Mu_2014}}
\bea
\begin{aligned}
	& \underset{\mc{X}}{\text{min}}&&||X_{[\text{round}(N/2)]}||_{*}~~~~ \text{s.t.}&& \mc{X}_{\Omega} = \mc{T}_{\Omega}.
\end{aligned}
\label{eq6_1}
\eea
by choosing the weights such that $\al_{k}=1$ if $k=\text{round}(N/2)$, otherwise $\al_{k}=0$. {\color{black}
Although the single matrix $X_{[\text{round}(N/2)]}$ is balanced and thus (\ref{eq6_1}) is an effective means for
minimizing $\mbox{rank}(X_{[\text{round}(N/2)]})$, it should be realized that  {\color{black}it only captures the
local correlation between $\text{round}(N/2)$ modes and other $\text{round}(N/2)$ modes.}}

The second problem is based on the factorization model
{\color{black}$X_{[k]} = UV$ for a matrix $X_{[k]}\in\mb{R}^{m\times n}$ of rank $r_k$,
where $U\in\mb{R}^{m\times r_k}$ and $V\in\mb{R}^{r_k\times n}$}. Instead of optimizing the
nuclear norm of the unfolding matrices $X_{[k]}$ {\color{black}as in (\ref{eq6})},
the Frobenius norm is minimized:
\bea
\begin{aligned}
	& \underset{U_k,V_k,\mc{X}}{\text{min}}&&\sum_{k=1}^{N-1}\frac{\al_{k}}{2}||U_kV_k-X_{[k]}||^{2}_{F} \\
	& \text{s.t.}&& \mc{X}_{\Omega} = \mc{T}_{\Omega},
\end{aligned}
\label{eq8}
\eea
where $U_{k}\in\mathbb{R}^{\prod_{j=1}^{k}I_{j}\times r_k}$ and $V_{k}\in\mathbb{R}^{r_k\times\prod_{j=k+1}^{N}I_{j}}$. This model is similar to the one proposed in \cite{Tan_2014,Xu} (which is an extension of the matrix completion model \cite{Wen_2012}) where the Tucker rank is employed.

\section{{\color{black}Proposed Algorithms}}\label{propalg}
{\color{black}This section is devoted to the algorithmic development for solutions of two optimization problems
(\ref{eq6}) and (\ref{eq8}).}
\subsection{SiLRTC-TT}
{\color{black}To address the problem (\ref{eq6}) we further convert it to the following problem:}
\bea
\begin{aligned}
	& \underset{\mc{X},M_k}{\text{min}}&&\sum_{k=1}^{N-1}\al_{k}||M_k||_{*} + \frac{\beta_k}{2}||X_{[k]}-M_k||^{2}_{F}\\
	& \text{s.t.}&& \mc{X}_{\Omega} = \mc{T}_{\Omega},
\end{aligned}
\label{eq7}
\eea
where $\bet_k$ are positive numbers. The central concept is based on the BCD method to alternatively optimize a group of variables while the other groups remain fixed. More specifically, the variables are divided into two main groups. The first one contains the unfolding matrices $M_1, M_2,\ldots, M_{N-1}$ and the other is tensor $\mc{X}$. Computing each matrix $M_k$ is related to solving the following optimization problem:
\bea
\begin{aligned}
	& \underset{M_k}{\text{min}}&&\al_{k}||M_k||_{*} + \frac{\beta_k}{2}||X_{[k]}-M_k||^{2}_{F},
\end{aligned}
\label{minrankX7}
\eea
with fixed $X_{[k]}$. The optimal solution for this problem has the closed form \cite{Ma_2009} which is determined by
\bea
M_k= \mathbf{D}_{\gamma_k}(X_{[k]}),
\eea
where $\gamma_k=\frac{\al_k}{\bet_k}$ and $ \mathbf{D}_{\gamma_k}(X_{[k]})$ denotes the thresholding SVD of $X_{[k]}$ \cite{Cai_2010}. Specifically, if the SVD of $X_{[k]} = U\lambda V^T$, its thresholding SVD is defined as:
\bea
\mathbf{D}_{\gamma_k}(X_{[k]}) = U\lambda_{\gamma_k}V^T,
\eea
where $\lambda_{\gamma_k} = diag(\max(\lambda_l-\gamma_k,0))$. After updating all the $M_k$ matrices, we turn into another block to compute the tensor $\mc{X}$ which elements are given by
\bea
x_{i_1\cdots i_N}= \left\{\begin{array}{ll}\Big (\frac{\sum_{k=1}^{N} \bet_k \text{fold}(M_k)}{\sum_{k=1}^{N}\bet_k}\Big)_{i_1\cdots i_N}&({i_1\cdots i_N})\notin\Omega\\
t_{i_1\cdots i_N}&({i_1\cdots i_N})\in\Omega\\
\end{array}\right.
\label{UpdateX}
\eea
The pseudocode of this algorithm is given in Algorithm \ref{Algorithm 1}. {\color{black}We call it
simple low-rank tensor completion via tensor train (SiLRTC-TT) as it is an enhancement of SiLRTC \cite{Ji2013}}.
The convergence condition is reached when the relative error between two successive tensors $\mc{X}$ is smaller than a threshold. The algorithm is guaranteed to be converged and gives rise to a global solution since the objective in (\ref{eq7}) is a convex and the nonsmooth term is separable. We can also apply this algorithm for the square model \cite{Mu_2014} by simply choosing the weights such that $\al_{k}=1$ if $k=\text{round}(N/2)$ otherwise $\al_{k}=0$. For this particular case, the algorithm is defined as SiLRTC-Square.
\begin{table}[!thb]
	\centering
	\caption{SiLRTC-TT}
	\label{Algorithm 1}	
	\begin{tabular}{*2l} 
		\hline
		~&~\\
		\tb{Input:} The observed data $\mc{T}\in\mathbb{R}^{I_1\times I_2\cdots\times I_{N}}$, index set $\Omega$.&\\
		\tb{Parameters:} $\al_k,\bet_k, k=1,\ldots,N-1$.&\\							
		~&~\\	
		\hline		
		~&~\\	
		1:~~\tb{Initialization:} $\mc{X}^{0}$, with $\mc{X}^{0}_{\Omega} = \mc{T}_{\Omega}$, $l=0$.& \\
		2:~~\tb{While not converged do:}&\\
		3:~~~~\tb{for} $k = 1$ \tb{to} $N-1$ \tb{do}&\\
		4:~~~~~~~Unfold the tensor $\mc{X}^{l}$ to get $X^{l}_{[k]}$&\\
		5:~~~~~~~$M^{l+1}_k = \mathbf{D}_{\frac{\al_k}{\bet_k}}(X^{l}_{[k]})$&\\	
		6:~~~~\tb{end for}&\\	
		7:~~ Update $\mc{X}^{l+1}$ from $M^{l+1}_k$ by (\ref{UpdateX})&\\
		8:~~\tb{End while}&\\
		\hline
		~&~\\
		\tb{Output:} The recovered tensor $\mc{X}$ as an approximation of $\mc{T}$&
	\end{tabular}
\end{table}
\subsection{TMac-TT}
To solve the problem given by (\ref{eq8}), {\color{black} following  TMac and TC-MLFM in \cite{Xu} and \cite{Tan_2014}},
we apply the BCD method to alternatively optimize different groups of variables. Specifically, we focus on the following problem:
\bea
\begin{aligned}
	& \underset{U_k,V_k,X_{[k]}}{\text{min}}&&||U_kV_k-X_{[k]}||^{2}_{F},
\end{aligned}
\label{minrankX7_1}
\eea
for $k=1,2,\ldots,N-1$.  This problem is convex when each variable $U_k,V_k$ and $X_{[k]}$ is modified while keeping the other two fixed.  To update each variable, perform the following steps:
\bea
U^{l+1}_{k} &=& X^{l}_{[k]}(V^{l}_{k})^{T}(V^{l}_{k}(V^{l}_{k})^{T})^{\dagger},\label{U1}\\
V^{l+1}_{k} &=&( (U^{l+1}_{k})^{T}U^{l+1}_{k})^{\dagger}(U^{l+1}_{k})^{T})X^{l}_{[k]}\\
X^{l+1}_{[k]}&=& U^{l+1}_{k}V^{l+1}_{k},
\label{X1}
\eea
where \textquotedblleft$^\dagger$\textquotedblright denotes the Moore-Penrose pseudoinverse. It was shown in \cite{Xu} that we can replace (\ref{U1}) by the following:
\bea
U^{l+1}_{k} &=& X^{l}_{[k]}(V^{l}_{k})^{T},\label{U2}
\eea
to avoid computing the Moore-Penrose pseudoinverse $(V^{l}_{k}(V^{l}_{k})^{T})^{\dagger}$. The rationale behind this is that we only need the product $U^{l+1}_{k}V^{l+1}_{k}$ to compute $X^{l+1}_{[k]}$ in (\ref{X1}), which is the same when either (\ref{U1}) or (\ref{U2}) is used. After updating $U^{l+1}_k,V^{l+1}_k$ and $X^{l+1}_{[k]}$ for all $k=1,2,\ldots,N-1$, we compute elements of the tensor $\mc{X}^{l+1}$ as follows:
\bea
x^{l+1}_{i_1\cdots}= \left\{\begin{array}{ll}\Big(\sum\limits_{k=1}^{N-1} \al_k \text{fold}(X^{l+1}_{[k]})\Big)_{i_1\cdots}&({i_1\cdots})\notin\Omega\\
	t_{i_1\cdots}&({i_1\cdots})\in\Omega
\end{array}\right.
\label{UpdateX2}
\eea
This algorithm is defined as tensor completion by parallel matrix factorization in the concept of tensor train (TMac-TT), and its pseudocode is summarized in Algorithm \ref{Algorithm 2}. {\color{black}  The essential advantage of this algorithm  is that
it avoids a lot of SVDs, and hence it can substantially save computational time.}

The algorithm can also be applied for the square model \cite{Mu_2014} by choosing the weights such that $\al_{k}=1$ if $k=\text{round}(N/2)$, otherwise $\al_{k}=0$. For this case, we define the algorithm TMac-Square.
\begin{table}[!thb]
	\centering
	\caption{TMac-TT}
	\label{Algorithm 2}	
	\begin{tabular}{*2l} 
		\hline
		~&~\\
		\tb{Input:} The observed data $\mc{T}\in\mathbb{R}^{I_1\times I_2\cdots\times I_{N}}$, index set $\Omega$.&\\
		\tb{Parameters:} $\al_i,r_i, i=1,\ldots,N-1$.&\\							
		~&~\\	
		\hline		
		~&~\\	
		1:~~\tb{Initialization:} $U^{0}, V^{0}, \mc{X}^{0}$, with $\mc{X}^{0}_{\Omega} = \mc{T}_{\Omega}$, $l=0$.& \\
		\tb{While not converged do:}&\\
		2:~~\tb{for} $k = 1$ \tb{to} $N-1$ \tb{do}&\\
		3:~~~~~~~Unfold the tensor $\mc{X}^{l}$ to get $X^{l}_{[k]}$&\\
		4:~~~~~~~$U^{l+1}_{i} = X^{l}_{[k]}(V^{l}_{k})^{T}$&\\
		5:~~~~~~~$V^{l+1}_{k} = ((U^{l+1}_{k})^{T}U^{l+1}_{k})^{\dagger}(U^{l+1}_{k})^{T}X^{l}_{[k]}$&\\	
		6:~~~~~~~$X^{l+1}_{[k]} = U^{l+1}_{k}V^{l+1}_{k}$&\\		
		7:~~\tb{end}&\\	
		8:~~Update the tensor $\mc{X}^{l+1}$ using (\ref{UpdateX2})\\
		\tb{End while}&\\
		\hline
		~&~\\
		\tb{Output:} The recovered tensor $\mc{X}$ as an approximation of $\mc{T}$&
	\end{tabular}
\end{table}

\subsection{Computational complexity of algorithms}
The computational complexity of the algorithms are given in Table \ref{Cost} to complete a tensor $\mc{X}\in\mathbb{R}^{I_1\times I_2\times\cdots\times I_N}$, where we assume that $I_1=I_2=\cdots= I_{N}=I$. The Tucker rank and TT rank are assumed to be equal, i.e. $r_1=r_2=\cdots= r_{N}=r$.
\captionsetup[table]{name={\bf Table}}
\setcounter{table}{0}
\begin{table}[!h]
	\caption{Computational complexity of algorithms for one iteration.}
	\label{Cost}
	\centering 
	\begin{tabular}{l l}
		Algorithm & Computational complexity\\
		\hline
		~\\
		SiLRTC  & $O(NI^{N+1})$\\
		SiLRTC-TT & $O(I^{3N/2}+I^{3N/2-1})$\\
		TMac & $O(3NI^Nr)$\\
		TMac-TT & $O(3(N-1)I^Nr)$\\
		\hline
	\end{tabular}
\end{table}

\section{{\color{black}Tensor augmentation} \label{sec3}}
In this section, we introduce \emph{ket augmentation} (KA) to represent a low-order tensor by a higher-order one, i.e.
 to  cast an $N$th-order tensor
$\mc{T}\in\mathbb{R}^{I_1\times I_2\times\cdots\times I_N}$ into a $K$th-order tensor
$\tilde{\mc{T}}\in\mathbb{R}^{J_1\times J_2\times\cdots\times J_K}$, where $K\geq N$ and $\prod_{l=1}^{N}I_l=\prod_{l=1}^{K}J_l$. A higher-order representation of the tensor offers some important advantages. For instance, the TT decomposition is more efficient for
the augmented tensor because the local structure of the data can be exploited effectively in terms of computational resources.
Actually, if the tensor is slightly correlated, its augmented tensor can be represented by a low-rank TT
\cite{Oseledets_2011,Lattore}.

The concept of KA was originally introduced in \cite{Lattore} for casting a grayscale image into \emph{real ket state} of a Hilbert space, which is simply a higher-order tensor, using an appropriate block structured addressing.

We define KA as a generalization of the original scheme to third-order tensors $\mc{T}\in\mathbb{R}^{I_1\times I_2\times I_{3}}$ that represent color images, where $I_1\times I_2=2^n\times 2^n$ ($n\geq 1\in\mathbb{Z}$) is the number of pixels in the image and $I_3 = 3$ is the number of colors (red, green and blue).
Let us start with an initial block, labeled as $i_1$, of $2\times 2$ pixels corresponding to a single color $j$ (assume that the color is indexed by $j$ where $j=1,2,3$ corresponding to red, green and blue colors, respectively). This block can be represented as
\bea
\mc{T}_{[2^1\times 2^1\times 1]} = \sum_{i_1=1}^{4}c_{i_1j}\tb{e}_{i_1},
\eea
where $c_{i_1j}$ is the pixel value corresponding to color $j$ and $\tb{e}_{i_1}$ is the orthonormal base which is defined as $\tb{e}_{1} = (1,0,0,0)$, $\tb{e}_{2} = (0,1,0,0)$, $\tb{e}_{3} = (0,0,1,0)$ and $\tb{e}_{4} = (0,0,0,1)$. The value $i_1=1,2,3$ and $4$ can be considered as labeling the up-left, up-right, down-left and down-right pixels, respectively. For all three colors, we have three blocks which are presented by
\bea
\mc{T}_{[2^1\times 2^1\times 3]} = \sum_{i_1=1}^{4}\sum_{j=1}^{3}c_{i_1j}\tb{e}_{i_1}\otimes \tb{u}_{j},
\label{ket1}
\eea
where $\tb{u}_{j}$ is also an orthonormal base which is defined as $\tb{u}_{1} = (1,0,0)$, $\tb{u}_{2} = (0,1,0)$, $\tb{u}_{3} = (0,0,1)$.
\begin{figure}[htpb]
	\centering
	\includegraphics[width=\columnwidth]{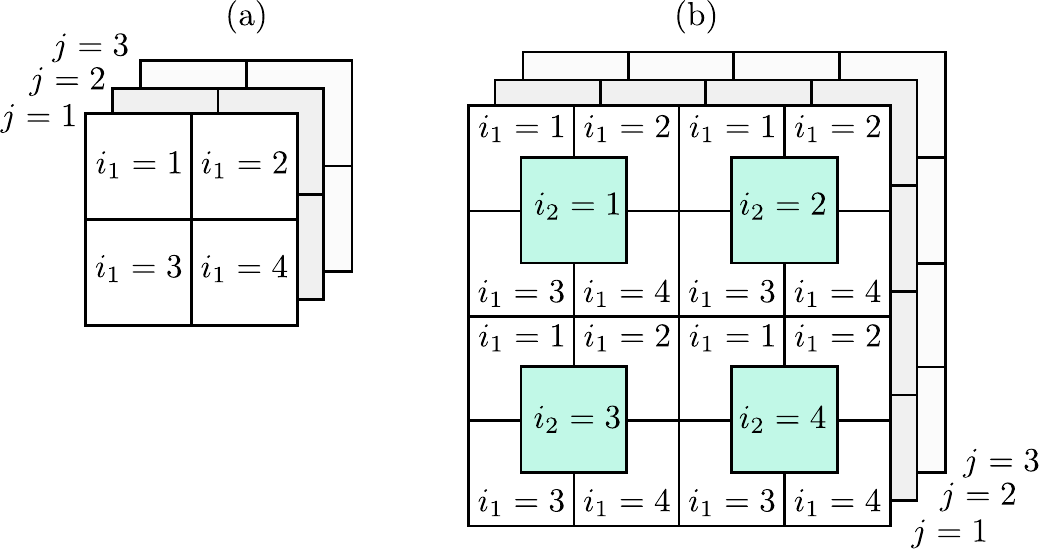}\\
	\caption{A structured block addressing procedure to cast an image into a higher-order tensor. (a) Example for an image of size $2\times 2\times 3$ represented by (\ref{ket1}). (b) Illustration for an image of size $2^2\times 2^2\times 3$ represented by (\ref{ket2}).}
	\label{fig10}
\end{figure}
We now consider a larger block labeled as $i_2$ make up of four inner sub-blocks for each color $j$ as shown in Fig.~\ref{fig10}. In total, the new block is represented by
\bea
\mc{T}_{[2^2\times 2^2\times 3]} = \sum_{i_2=1}^{4}\sum_{i_1=1}^{4}\sum_{j=1}^{3}c_{i_2i_1j}\tb{e}_{i_2}\otimes\tb{e}_{i_1}\otimes \tb{u}_{j}.
\label{ket2}
\eea
Generally, this block structure can be extended to a size of $2^n\times 2^n\times 3$ after several steps until it can present all the values of pixels in the image. Finally, the image can be cast into an $(n+1)$th-order tensor $\mc{C}\in\mb{R}^{4\times 4\times\cdots\times 4\times 3}$ containing all the pixel values as follows,
\bea
\mc{T}_{[2^n\times 2^n\times 3]} = \sum_{i_n,\ldots,i_1=1}^{4}\sum_{j=1}^{3}c_{i_n\cdots i_1j}\tb{e}_{i_n}\otimes\cdots\otimes\tb{e}_{i_1}\otimes \tb{u}_{j}.
\eea
This presentation is suitable for image processing as it not only preserves the pixels values, but also rearranges them in a higher-order tensor such that the richness of textures in the image can be studied via the correlation between modes of the tensor \cite{Lattore}. Therefore, due to the flexibility of the TT-rank, our proposed algorithms would ideally take advantage of KA.

\section{Simulations \label{sec4}}
Extensive experiments are conducted with synthetic data, color images and videos. The proposed algorithms are benchmarked against TMac \cite{Xu}, TMac-Square, SiLRTC \cite{Ji2013}, SiLRTC-Square \cite{Mu_2014} and state-of-the-art tensor completion methods FBCP \cite{7010937} and STDC \cite{6587455}\footnote{Applicable only for tensors of order $N=3$.}. Additionally, we also benchmark the TT-rank based optimization algorithm, ALS \cite{Holtz_2012,Grasedyck_2015}.

The simulations for the algorithms are tested with respect to different missing ratios ($mr$) of the test data, with $mr$ defined as
\bea
mr = \frac{p}{\prod_{k=1}^{N}I_k},
\eea
where $p$ is the number of missing entries, which is chosen randomly from a tensor $\mc{T}$ based on a uniform distribution.

To measure performance of a LRTC algorithm, the relative square error (RSE) between the approximately recovered tensor $\mc{X}$ and the original one $\mc{T}$ is used, which is defined as,
\bea
RSE = ||\mc{X}-\mc{T}||_F/||\mc{T}||_F.
\label{RSE}
\eea

The convergence criterion of our proposed algorithms is defined by computing the relative error of the tensor $\mc{X}$ between two successive iterations as follows:
\bea
\epsilon = \frac{||\mc{X}^{l+1} - \mc{X}^{l}||_{F}}{||\mc{T}||_{F}}\leq tol,
\eea
where $tol=10^{-4}$ and the maximum number of iterations $maxiter = 1000$.
These simulations are implemented under a Matlab environment.

\subsection{Initial parameters}
In the experiments there are three parameters that must be initialized: the weighting parameters $\al$ and $\bet$, and the initial TT ranks ($r_i,i=1,\ldots,N-1$) for TMac, TMac-TT and TMac-Square. Firstly, the weights $\al_k$ are defined as follows:
\bea
\al_k &=& \frac{\delta_{k}}{\sum_{k=1}^{N-1}\delta_{k}}~~\text{with}~~\delta_{k} = \min(\prod_{l=1}^{k}I_l,\prod_{l=k+1}^{N}I_l),~~~
\eea
where $k=1,\ldots,N-1$. In this way, we assign the large weights to the more balanced matrices. The positive parameters are chosen by $\bet_k = f\al_k$, where $f$ is empirically chosen from one of the following values in $[0.01, 0.05, 0.1, 0.5, 1]$ in such a way that the algorithm performs the best. Similarly, for SiLRTC and TMac, the weights are chosen as follows:
\bea
\al_k &=& \frac{I_{k}}{\sum_{k=1}^{N}I_{k}},
\eea
where  $k=1,\ldots,N$. The positive parameters are chosen such that $\bet_k = f\al_k$, where $f$ is empirically chosen from one of the following values in $[0.01, 0.05, 0.1, 0.5, 1]$ which gives the best performance.

To obtain the initial TT ranks for TMac, TMac-TT and TMac-Square, each rank $r_i$ is bounded by keeping only the singular values that satisfy the following inequality:
\beq
\frac{\lambda^{[i]}_j}{\lambda^{[i]}_1}>th,
\eeq
with $j = 1,\ldots,r_i$, threshold $th$, and $\{\lambda^{[i]}_j\}$ is assumed to be in descending order. This condition is chosen such that the matricizations with low-rank (small correlation) will have more singular values truncated. We also choose $th$ empirically based on the algorithms performance.

It is important to highlight that these initial parameters can affect the performance of the proposed algorithms. Consequently, the proposed algorithms performance may not necessarily be optimal and future work will need to be considered in determining the optimal TT rank and weights via automatic \cite{7010937} and/or adaptive methods \cite{Xu}.

\subsection{Synthetic data completion}
We firstly perform the simulation on two different types of low-rank tensors which are generated synthetically in such a way that the Tucker and TT rank are known in advance.
\subsubsection{Completion of low TT rank tensor}
The $N$th-order tensors $\mc{T}\in\mathbb{R}^{I_1\times I_2\cdots\times I_{N}}$ of TT rank $(r_1,r_2,\ldots, r_{N-1})$ are generated such that its elements is represented by a TT format \cite{Oseledets_2011}. Specifically, its elements is $t_{i_{1}i_{2}\ldots i_{N}} = A^{[1]}_{i_1}A^{[2]}_{i_2}\cdots A^{[N]}_{i_N}$, where $A^{[1]}\in\mathbb{R}^{I_1\times r_{1}}$, $A^{[N]}\in\mathbb{R}^{r_{N}\times I_N}$ and $\mc{A}^{[k]}\in\mathbb{R}^{r_{k-1}\times I_k\times r_{k}}$ with $k = 2,\ldots,N-1$ are generated randomly with respect to the standard Gaussian distribution $\mc{N}(0,1)$. For simplicity, in this paper we set all components of the TT rank the same and so does the dimension of each mode, i.e. $r_1=r_2=\cdots= r_{N-1}=r$ and $I_1=I_2=\cdots= I_{N}=I$.

\begin{figure}[htpb]
	\centering
	\includegraphics[width=\columnwidth]{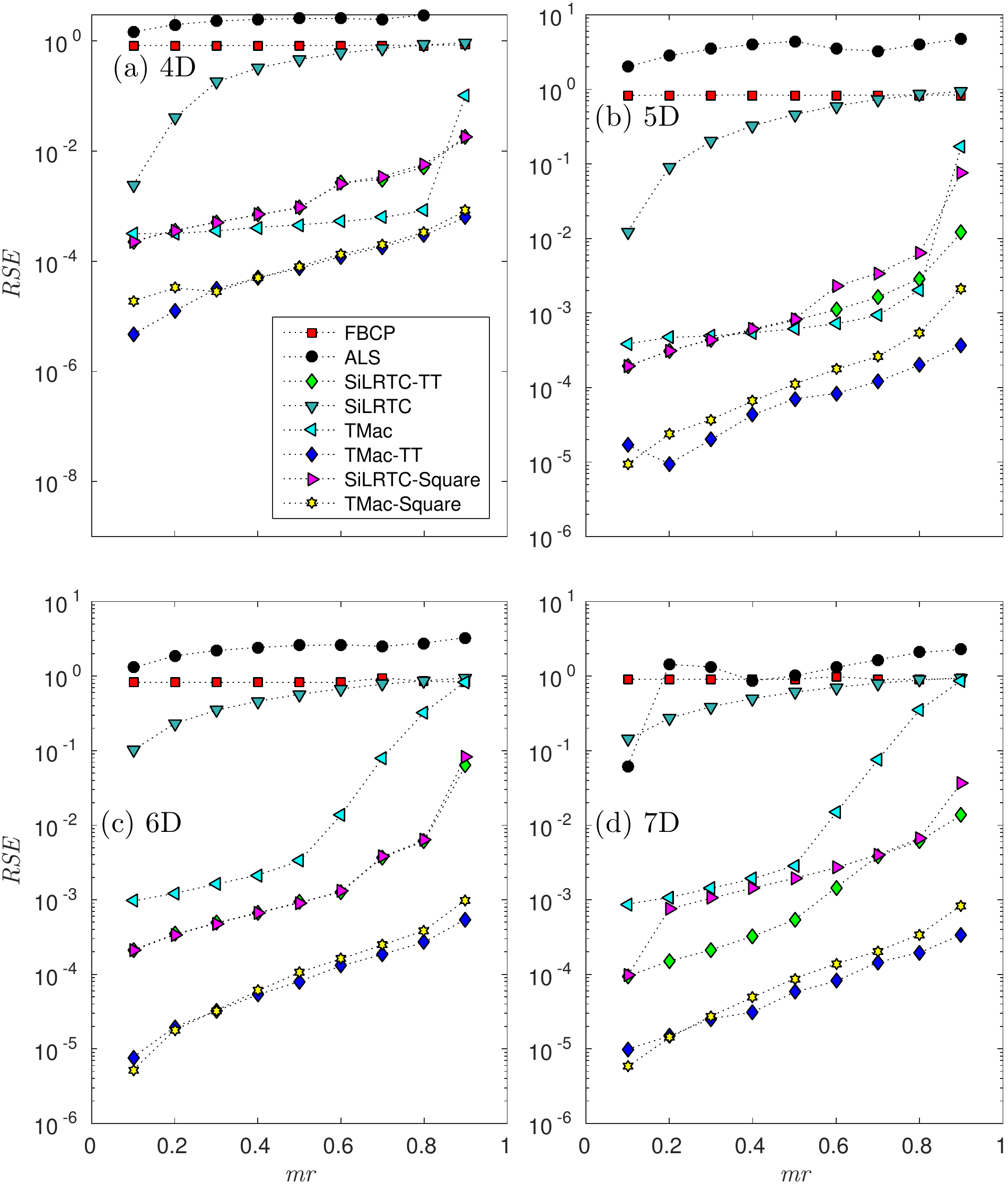}\\
	\caption{The RSE comparison when applying different LRTC algorithms to synthetic random tensors of low TT rank. Simulation results are shown for different tensor dimensions, 4D, 5D, 6D and 7D.}
	\label{fig1}
\end{figure}

The plots of RSE with respect to $mr$ are shown in the Figure.~\ref{fig1} for tensors of different sizes, $40\times 40\times 40\times 40$ (4D), $20\times 20\times 20\times 20\times 20$ (5D), $10\times 10\times 10\times 10\times 10\times 10$ (6D) and $10\times 10\times 10\times 10\times 10\times 10\times 10$ (7D) and the corresponding TT rank tuples are $(10, 10, 10)$ (4D), $(5, 5, 5,5)$ (5D), $(4,4,4,4,4)$ (6D) and $(4,4,4,4,4,4)$ (7D). From the plots we can see that TMac-TT shows best performance in most cases. Particularly, TMac-TT can recover the tensor successfully despite the high missing ratios, where in most cases with high missing ratios, e.g. $mr = 0.9$, it can recover the tensor with $RSE\approx 10^{-4}$. More importantly, the proposed algorithms SiLRTC-TT and TMac-TT often performs better than their corresponding counterparts, i.e. SiLRTC and TMac in most cases. FBCP and ALS have the worst results with random synthetic data, so for the remaining synthetic data experiments, only SiLRTC, SiLRTC-Square, SiLRTC-TT, TMac, TMac-Square and TMac-TT are compared.
\begin{figure}[htpb]
	\centering
	\includegraphics[scale = 0.5]{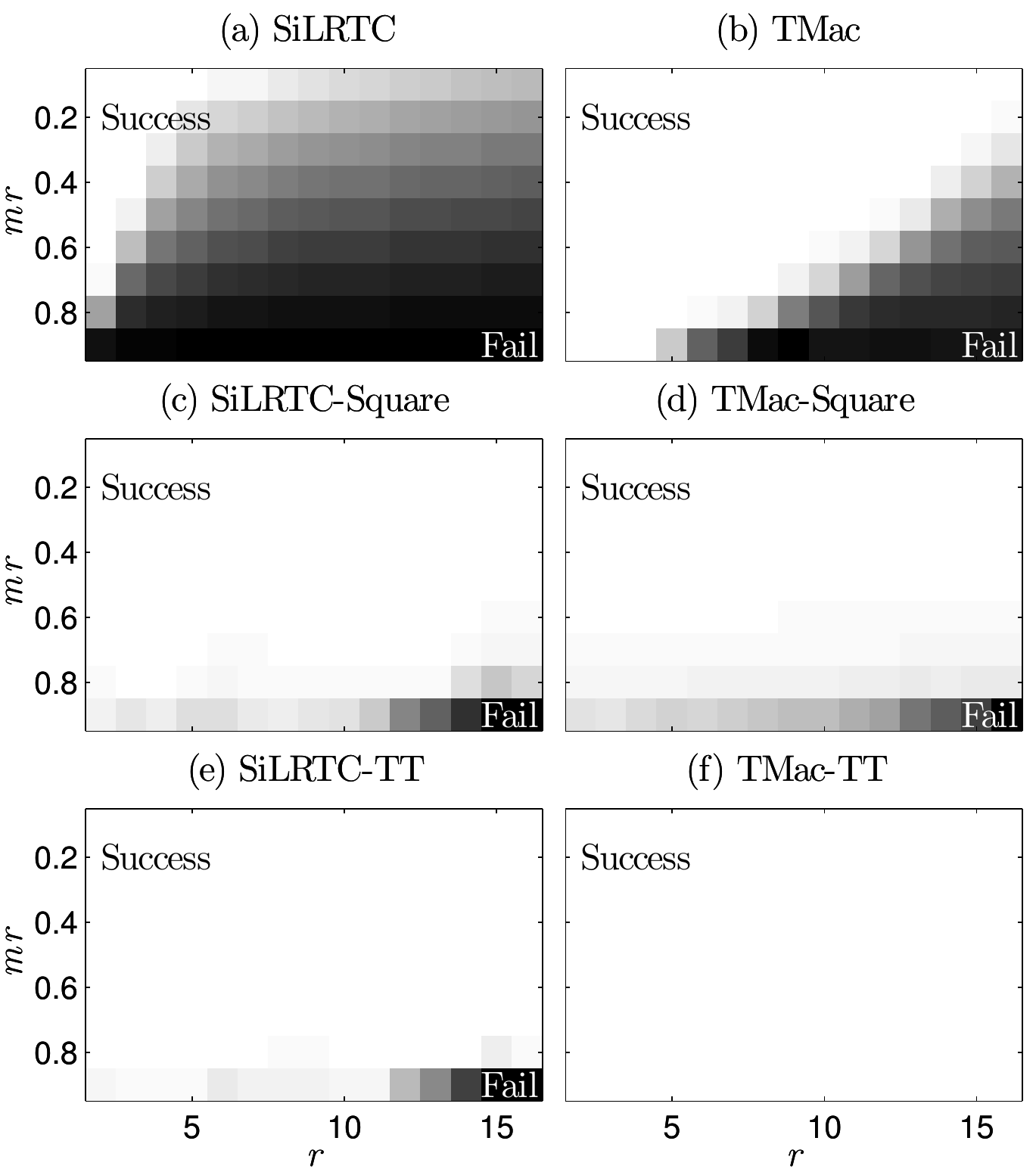}\\
	\caption{Phase diagrams for low TT rank tensor completion when applying different algorithms to a 5D tensor.}
	\label{fig2}
\end{figure}
\begin{figure*}
\begin{subfigure}{0.22\textwidth}
	\centering
	\includegraphics[scale = 0.4]{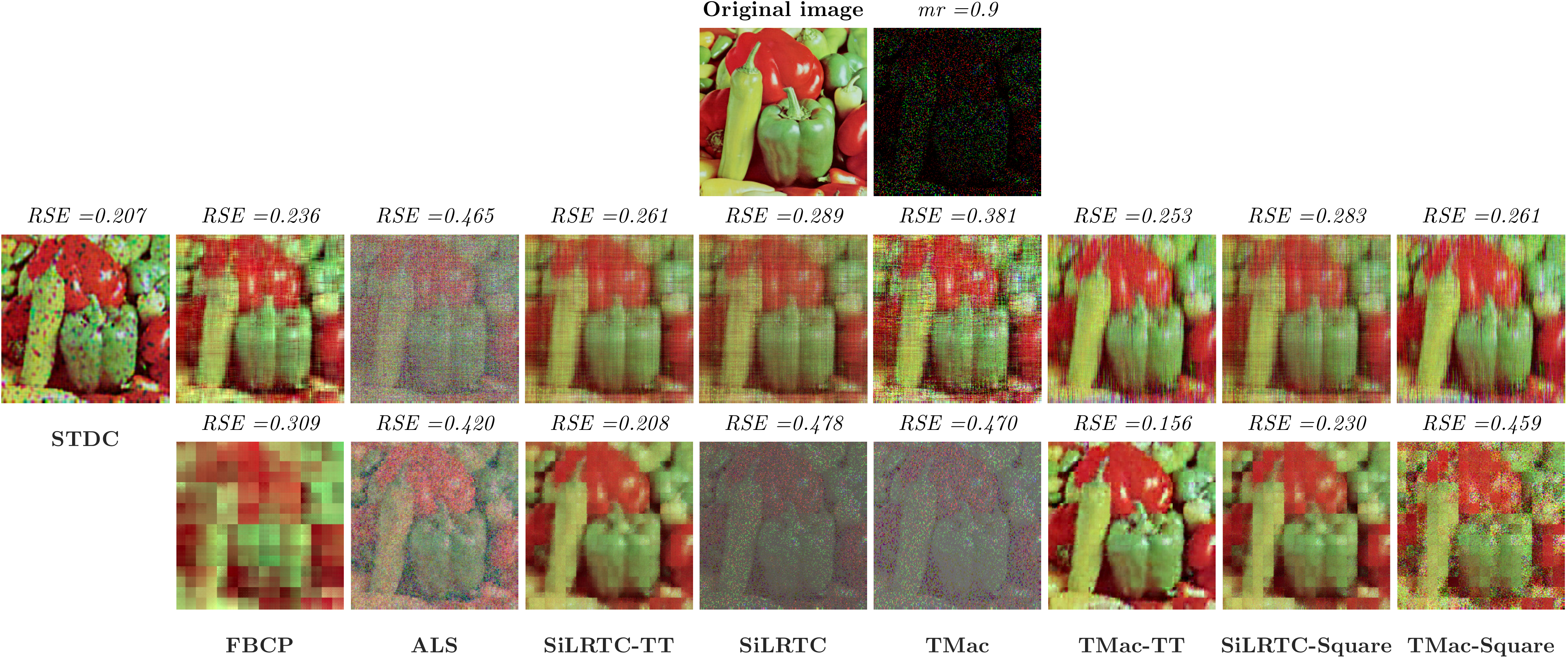}\\
\end{subfigure}\\
\caption{Recover the Peppers image with $90\%$ of missing entries using different algorithms. Top row from left to right: the original image and its copy with $90\%$ of missing entries. Second and third rows represent the recovery results of third-order (no order augmentation) and ninth-order tensors (KA augmentation), using different algorithms: STDC (only on second row), FBCP, ALS, SiLRTC-TT, SiLRTC, TMac, TMac-TT, SiLRTC-Square and TMac-Square from the left to the right, respectively.\newline}\label{fig6}
\end{figure*}

For a better comparison on the performance of different LRTC algorithms, we present the phase diagrams using the grayscale color to estimate how successfully a tensor can be recovered for a range of different TT rank and missing ratios. If $RSE\leq\epsilon$ where $\epsilon$ is a small threshold, we say that the tensor is recovered successfully and is represented by a white block in the phase diagram. Otherwise, if $RSE>\epsilon$, the tensor is recovered partially with a relative error and the block color is gray. Especially the recovery is completely failed if $RSE=1$. Concretely, we show in Fig.~\ref{fig2} the phase diagrams for different algorithms applied to complete a 5D tensor of size $20\times 20\times 20\times 20\times 20$ where the TT rank $r$ varies from 2 to 16 and $\epsilon=10^{-2}$. We can see that our LRTC algorithms outperform the others. Especially, TMac-TT always recovers successfully the tensor with any TT rank and missing ratio.

\begin{figure}[htpb]
	\centering
	\includegraphics[scale = 0.5]{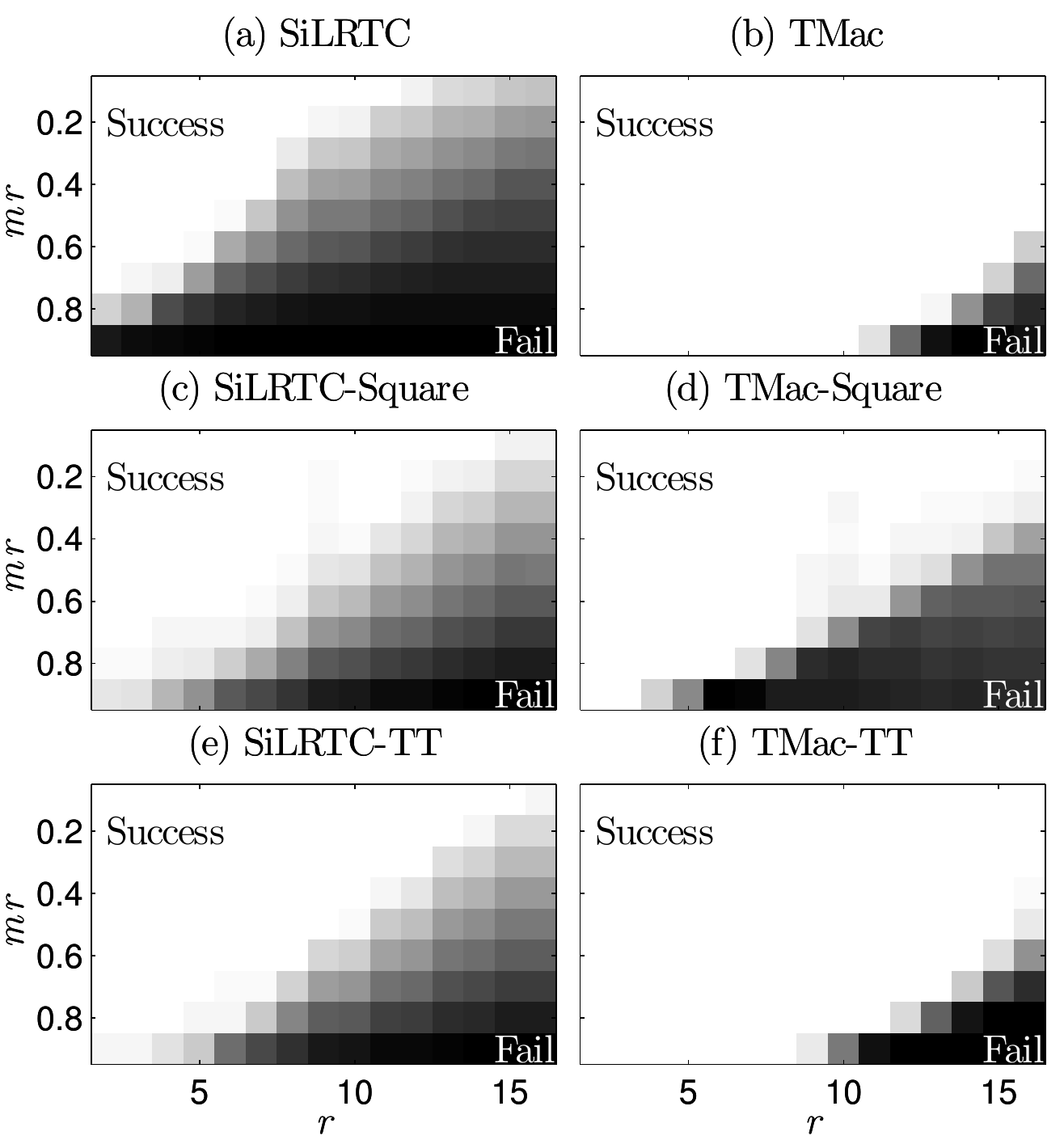}\\
	\caption{Phase diagrams for low Tucker rank tensor completion when applying different algorithms to a 5D tensor.}
	\label{fig4}
\end{figure}

\subsubsection{Completion of low Tucker rank tensor}
Let us now apply our proposed algorithms to synthetic random tensors of low Tucker rank. The $N$th-order tensor $\mc{T}\in\mathbb{R}^{I_1\times I_2\cdots\times I_{N}}$ of Tucker rank $(r_1,r_2,\ldots, r_N)$ is constructed by $\mc{T} = \mc{G}\times_1A^{(1)}\times_2A^{(2)}\cdots\times_NA^{(N)}$, where the core tensor $\mc{G}\in\mathbb{R}^{r_1\times r_2\cdots\times r_{N}}$ and the factor matrices $A^{(k)}\in\mathbb{R}^{r_k\times I_{k}},k=1,\ldots,N$ are generated randomly by using the standard Gaussian distribution $\mc{N}(0,1)$. Here, we choose $r_1=r_2=\cdots= r_{N}=r$ and $I_1=I_2=\cdots= I_{N}=I$ for simplicity. To compare the performance between the algorithms, we show in the Fig.~\ref{fig4} the phase diagrams for different algorithms applied to complete a 5D tensor of size $20\times 20\times 20\times 20\times 20$ where the Tucker rank $r$ varies from 2 to 16 and $\epsilon=10^{-2}$. We can see that both TMac and TMac-TT perform much better than the others. Besides, SiLRTC-TT shows better performance when compared to SiLRTC and SiLRTC-Square. Similarly, TMac-TT is better than its particular case TMac-Square.

In summary, we can see that although the tensors are generated synthetically to have low Tucker ranks, the proposed algorithms are still capable of producing results which are as good as those obtained by the Tucker-based algorithms.

The synthetic data experiments were performed to initially test the proposed algorithms. In order to have a better comparison between the algorithms we benchmark the methods against real world data such as color images and videos, where the ranks of the tensors are not known in advance. These will be seen in the subsequent subsections.

\begin{figure*}
\begin{subfigure}{0.22\textwidth}
	\centering
	\includegraphics[scale = 0.4]{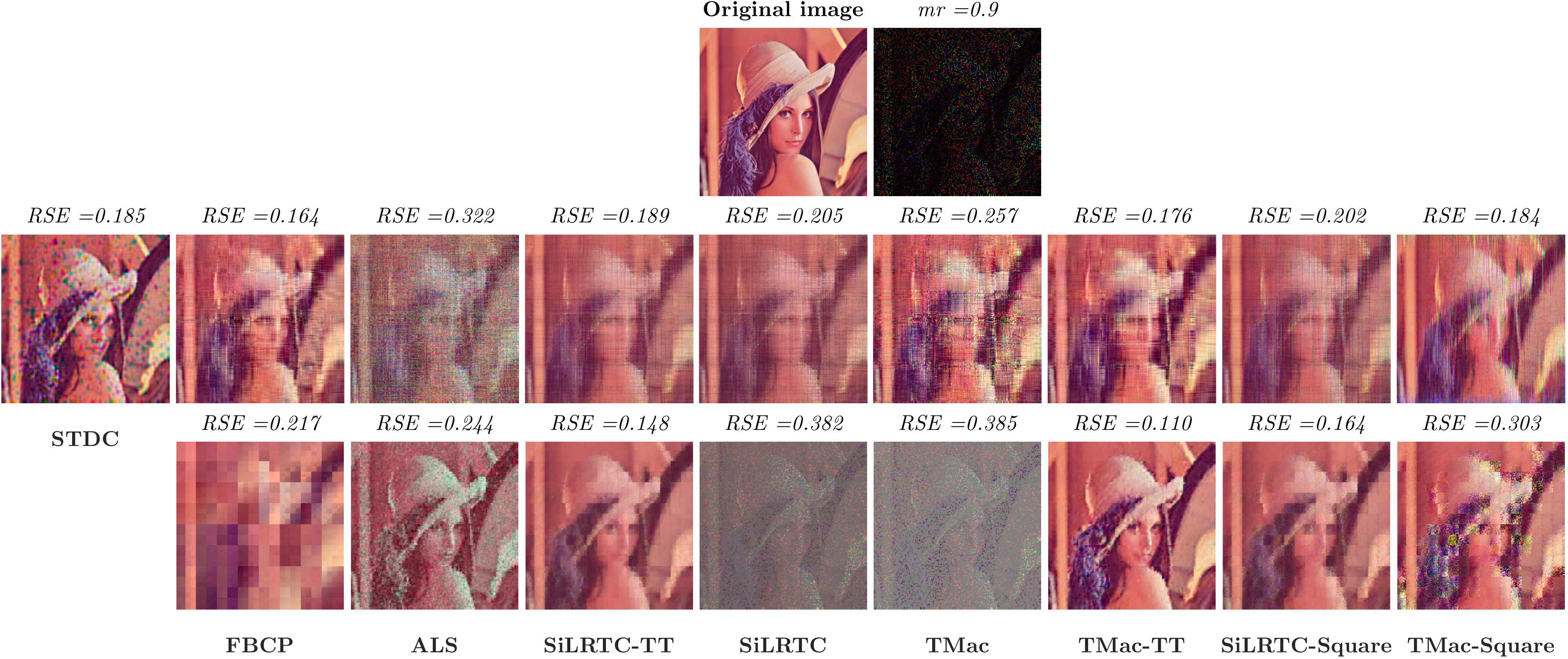}\\
\end{subfigure}\\
\caption{Recover the Lena image with $90\%$ of missing entries using different algorithms. Top row from left to right: the original image and its copy with $90\%$ of missing entries. Second and third rows represent the recovery results of third-order (no order augmentation) and ninth-order tensors (KA augmentation), using different algorithms: STDC (only on second row), FBCP, ALS, SiLRTC-TT, SiLRTC, TMac, TMac-TT, SiLRTC-Square and TMac-Square from the left to the right, respectively.}\label{fig8}
\end{figure*}
\begin{figure*}[htpb]
	\centering
	\includegraphics[scale = 0.4]{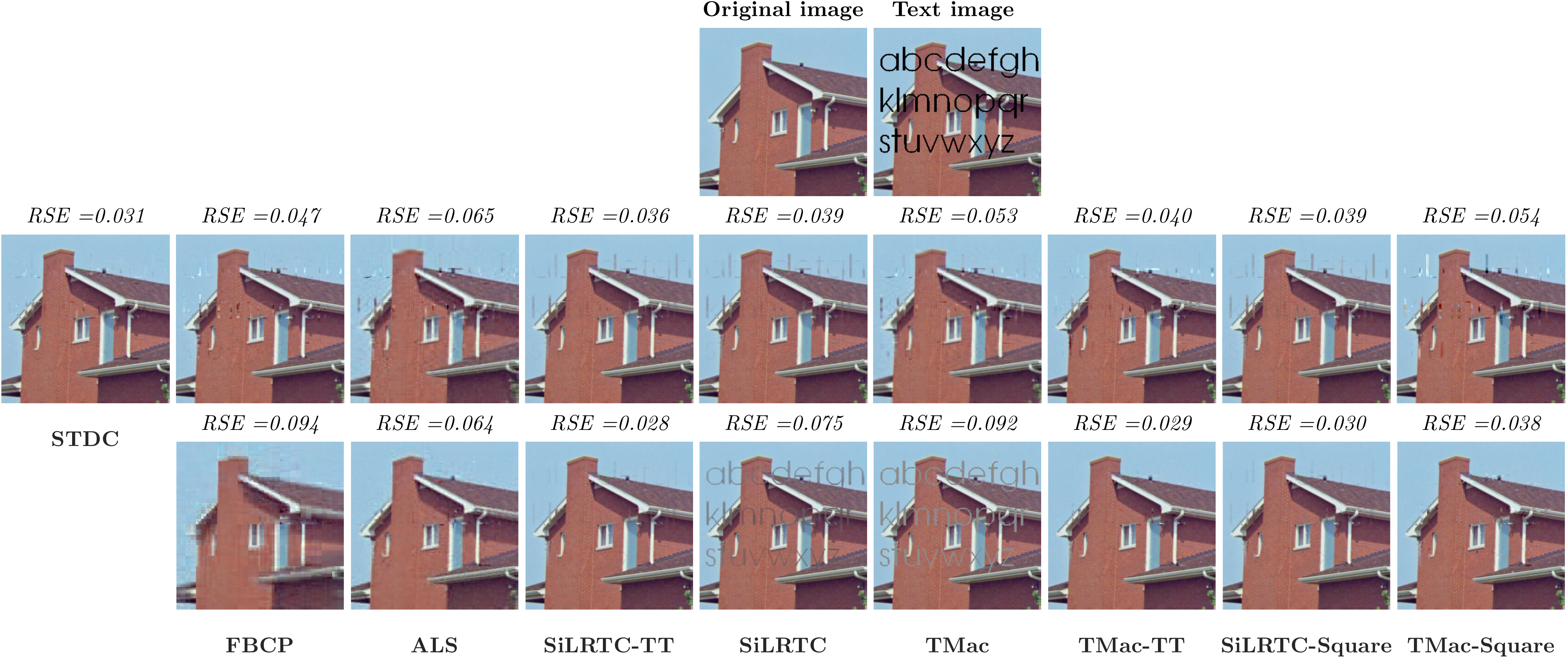}\\
	\caption{Recover the House image with missing entries described by the white letters using different algorithms. Top row from left to right: the original image and its copy with white letters. Second and third rows represent the recovery results of third-order (no order augmentation) and ninth-order tensors (KA augmentation), using different algorithms: STDC (only on second row), FBCP, ALS, SiLRTC-TT, SiLRTC, TMac, TMac-TT, SiLRTC-Square and TMac-Square from the left to the right, respectively.}
	\label{fig9}
\end{figure*}
		\begin{figure*}[h]
		\centering
		\begin{subfigure}{0.22\textwidth}
			\includegraphics[scale=0.13]{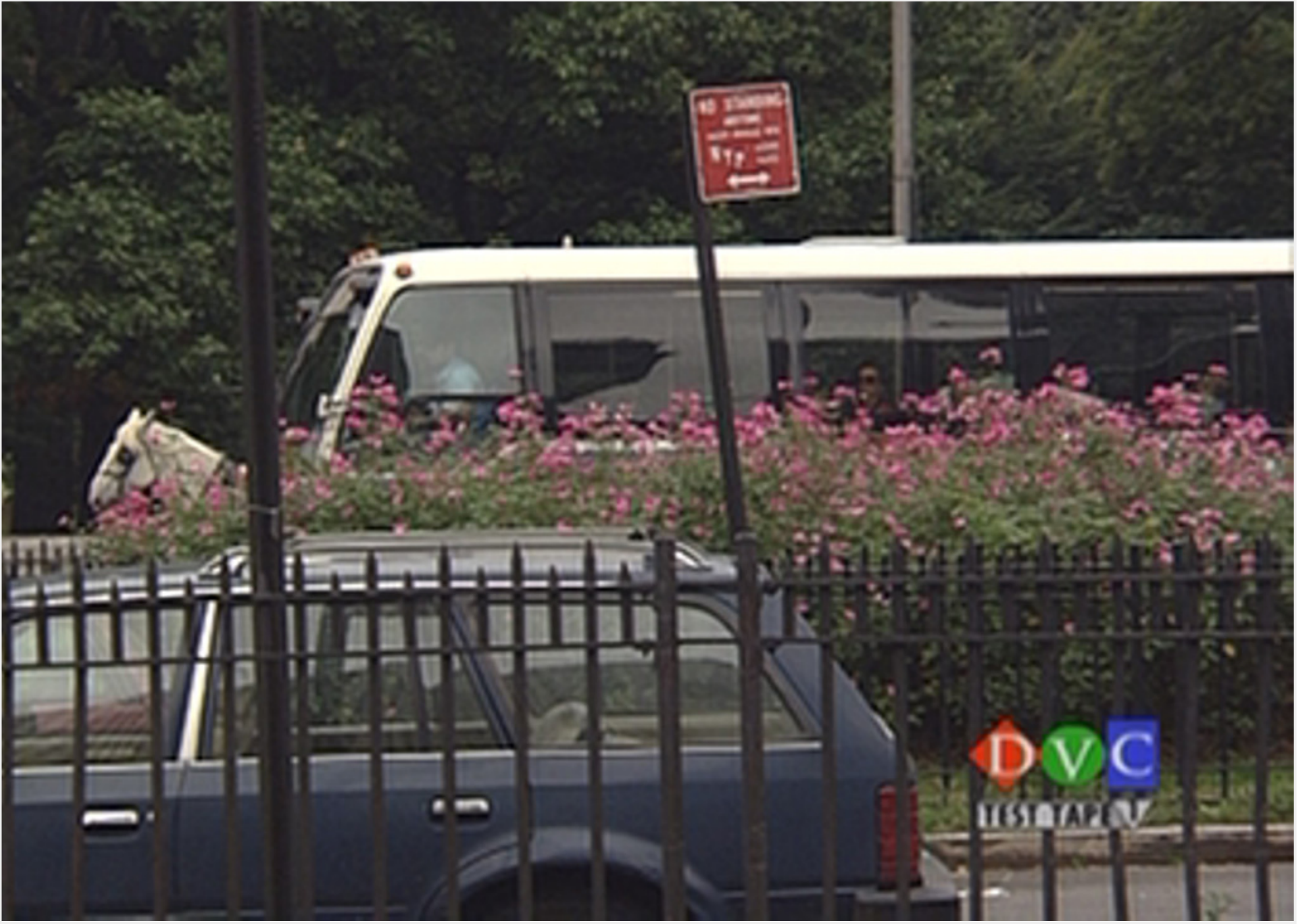}
			\caption{Bus frame 1}
			\label{fig:a}
		\end{subfigure}
		\centering
		\begin{subfigure}{0.22\textwidth}
			\includegraphics[scale=0.13]{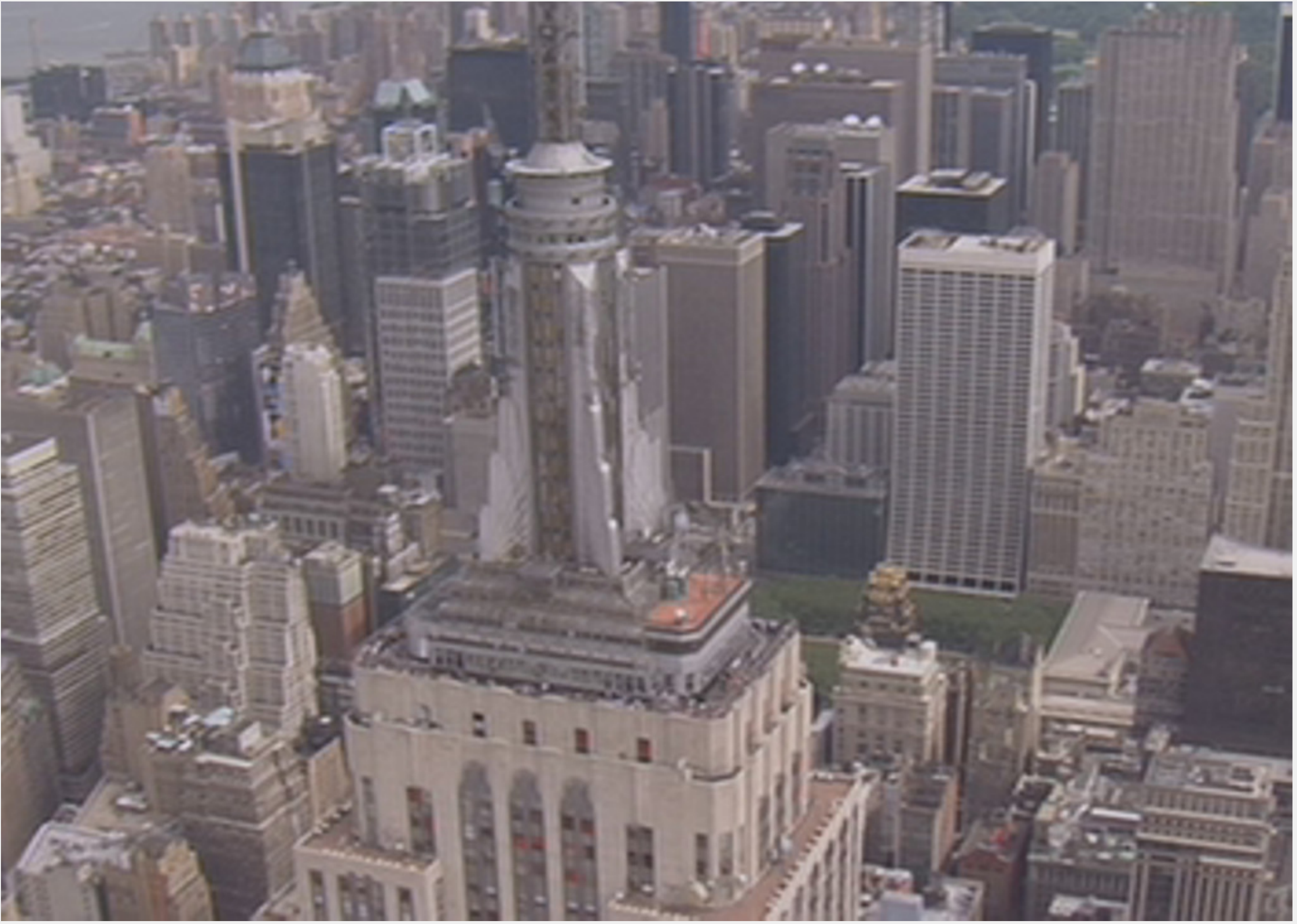}
			\caption{NYC frame 1}
			\label{fig:b}
		\end{subfigure}	
		\centering
		\begin{subfigure}{0.22\textwidth}
			\includegraphics[scale=0.13]{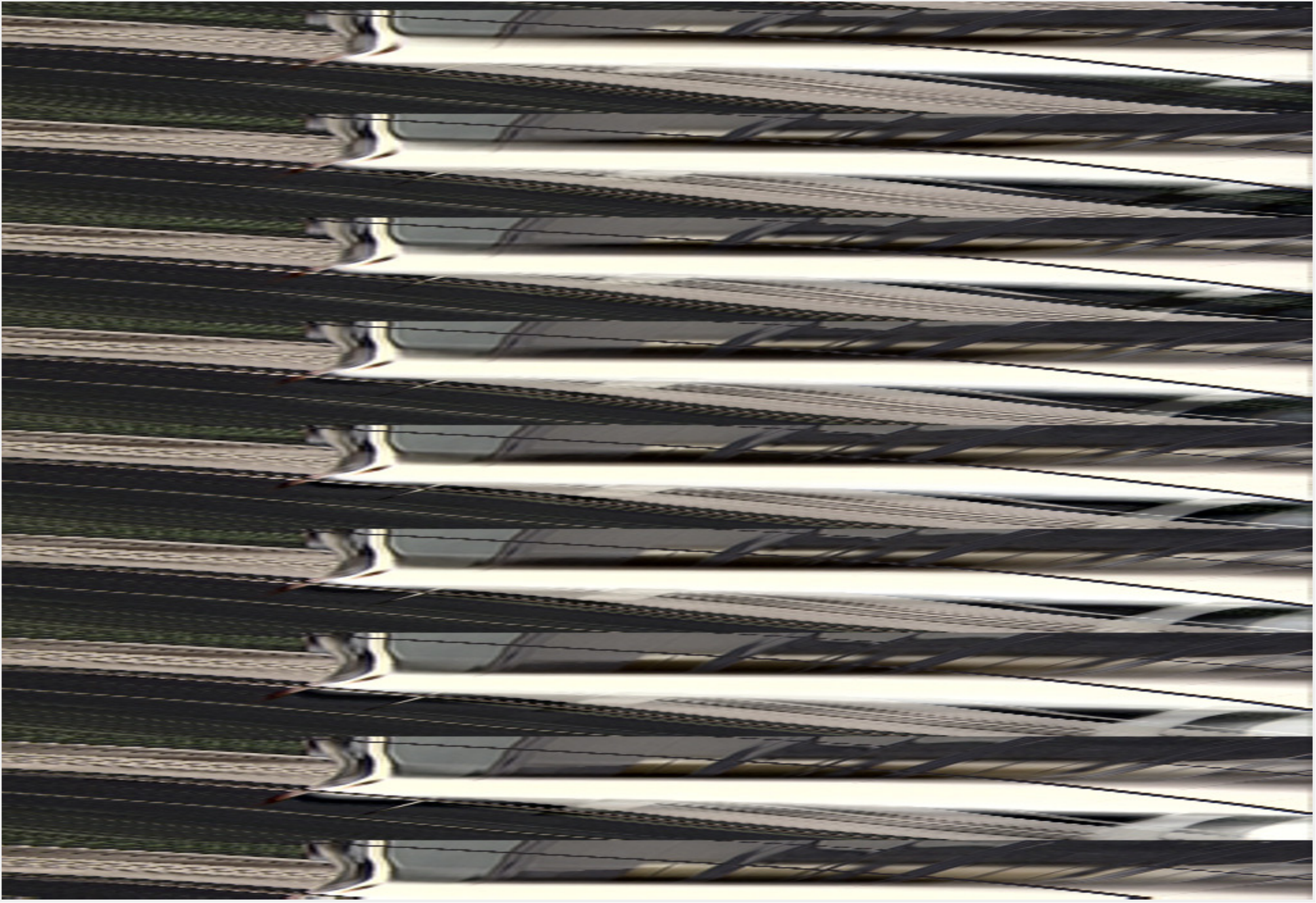}
			\caption{Bus combined (20000:20700)}
			\label{fig:c}
		\end{subfigure}	
		\centering	
		\begin{subfigure}{0.22\textwidth}
			\includegraphics[scale=0.13]{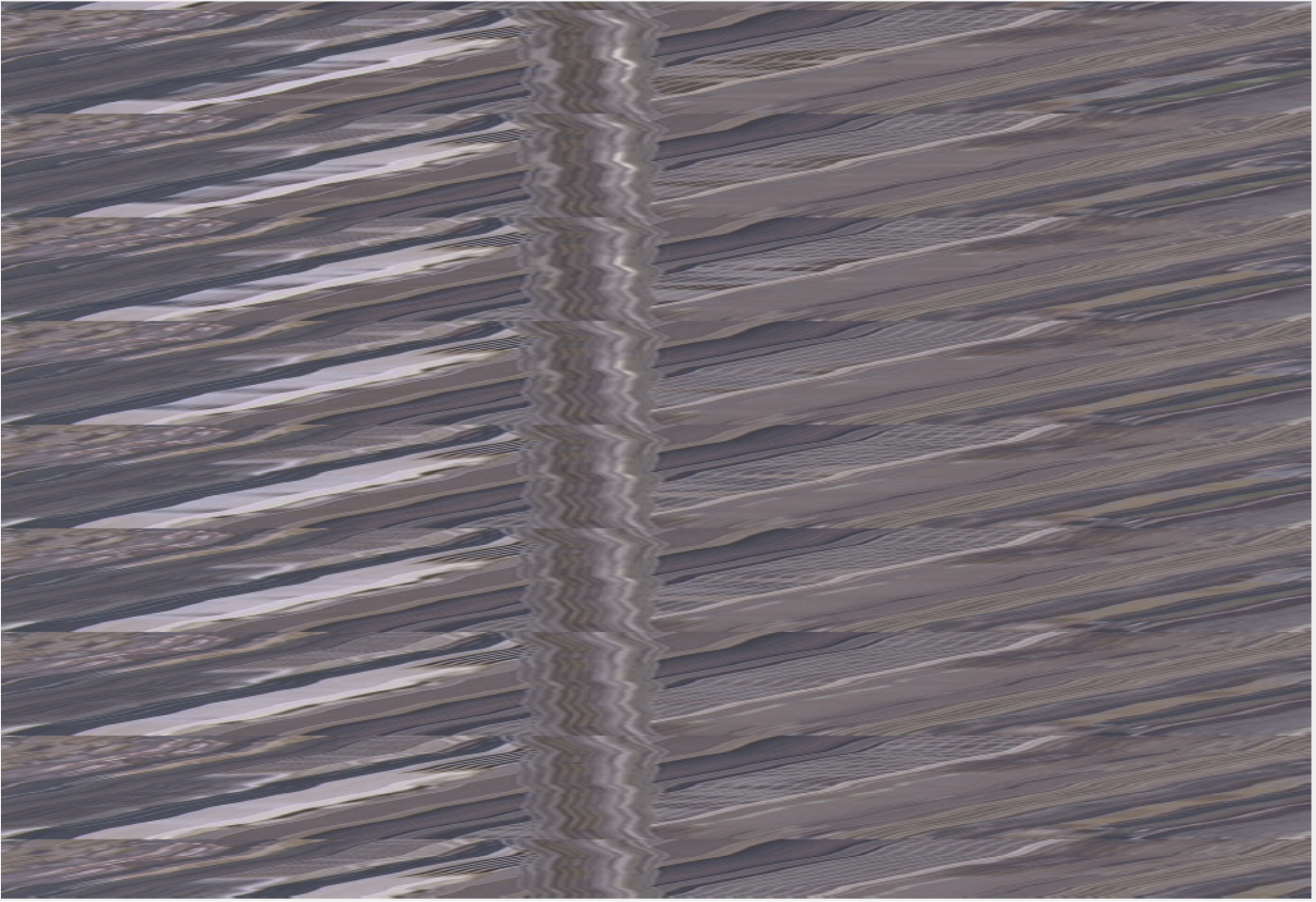}
			\caption{NYC combined (20000:20700)}
			\label{fig:d}
		\end{subfigure}					
		\caption{The first frames of the bus and NYC videos are shown (a) and (b), respectively. In (c) and (d), the third-order VST for bus and NYC are shown for the range 20000:20700 in the $combined\ row$ mode, respectively.}\label{dctfigs}
		\end{figure*}
\subsection{Image completion}
The color images known as Peppers, Lena and House are employed to test the algorithms. All the images are initially represented by third-order tensors which have same sizes of $256\times 256\times 3$.

\begin{figure}[htpb]
	\centering
	\includegraphics[width=\columnwidth]{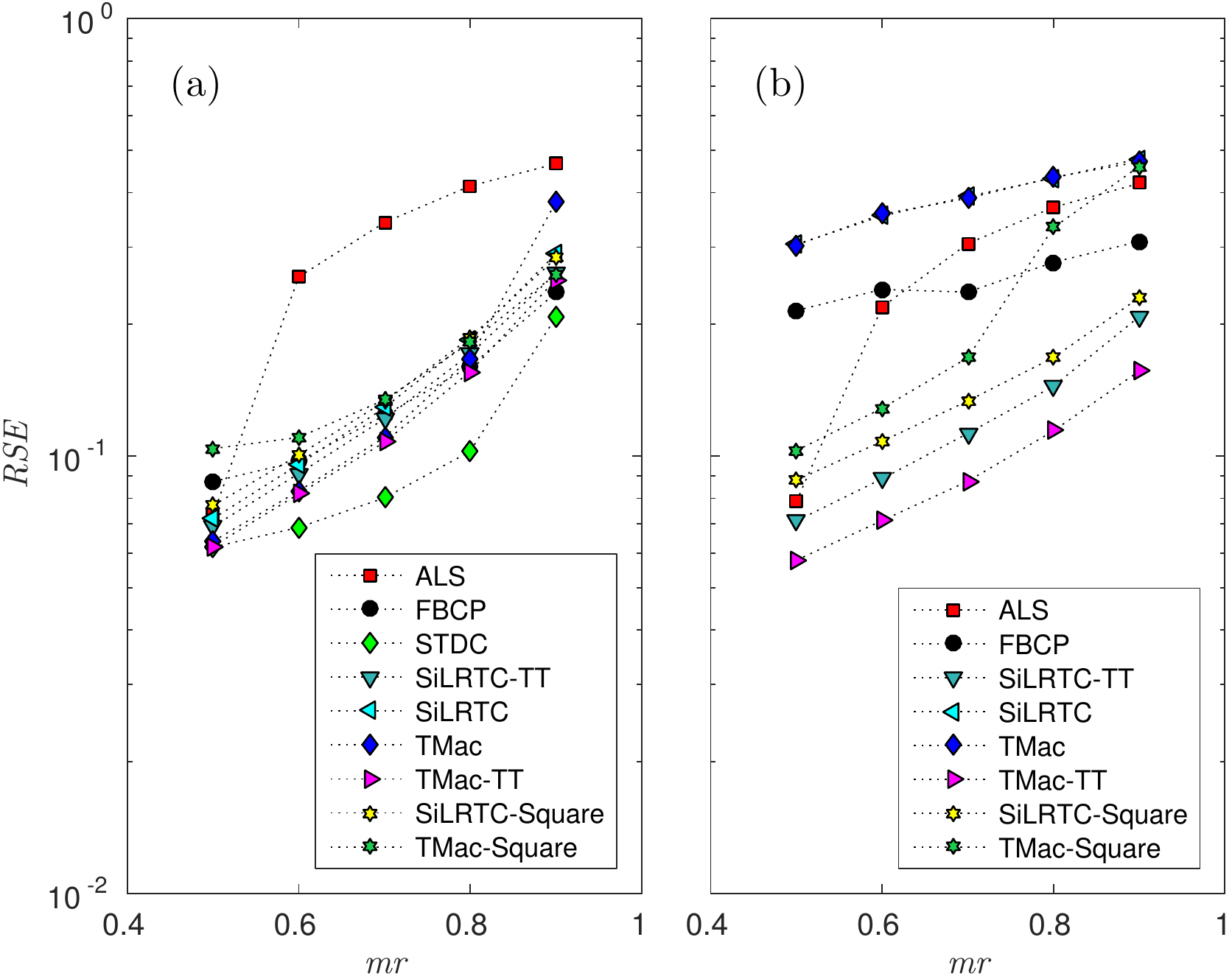}\\
	\caption{Performance comparison between different tensor completion algorithms based on the RSE vs the missing rate when applied to the Peppers image. (a) Original tensor (no order augmentation). (b) Augmented tensor using KA scheme.}
	\label{fig5}
\end{figure}
Note that when completing the third-order tensors, we do not expect the proposed methods to prevail against the conventional ones due to the fact that the TT rank of the tensor is a special case of the Tucker rank. Thus, performance of the algorithms should be mutually comparable. However, for the purpose of comparing the performance between different algorithms for real data (images) represented in terms of higher-order tensors, we apply the tensor augmentation scheme KA mentioned above to reshape third-order tensors to  higher-order ones without changing the number of entries in the tensor. Specifically, we start our simulation by casting a third-order tensor $\mc{T}\in\mathbb{R}^{256\times 256\times 3}$ into a ninth-order $\tilde{\mc{T}}\in\mathbb{R}^{4\times  4\times 4\times 4\times 4\times  4\times 4\times 4\times 3}$ and then apply the tensor completion algorithms to impute its missing entries. We perform the simulation for the Peppers and Lena images where missing entries of each image are chosen randomly according to a uniform distribution, the missing ratio $mr$ varies from 0.5 to 0.9.

In Fig.~\ref{fig5}, performance of the algorithms on completing the Peppers image is shown. When the image is represented by a third-order tensor, the STDC algorithm performs very well against all methods, with the ALS algorithm performing poorly, and the remaining algorithms having similar performance. However, for the case of the ninth-order tensors, the performance of the algorithms are rigorously distinguished. Specifically, our proposed algorithms (especially TMac-TT) prevails against all other methods, and this is demonstrated in Fig.~\ref{fig6} for $mr=0.9$. This shows that our proposed algorithms give really good results in the case of augmented tensors. Furthermore, using the KA scheme to increase the tensor order, SiLRTC-TT and TMac-TT are \textit{at least} comparable to STDC, with TMac-TT having the lowest $RSE$ for $mr=0.9$. More precisely, TMac-TT gives the best result of $RSE\approx0.156$ when using the KA scheme.

\begin{figure}[htpb]
	\centering
	\includegraphics[width=\columnwidth]{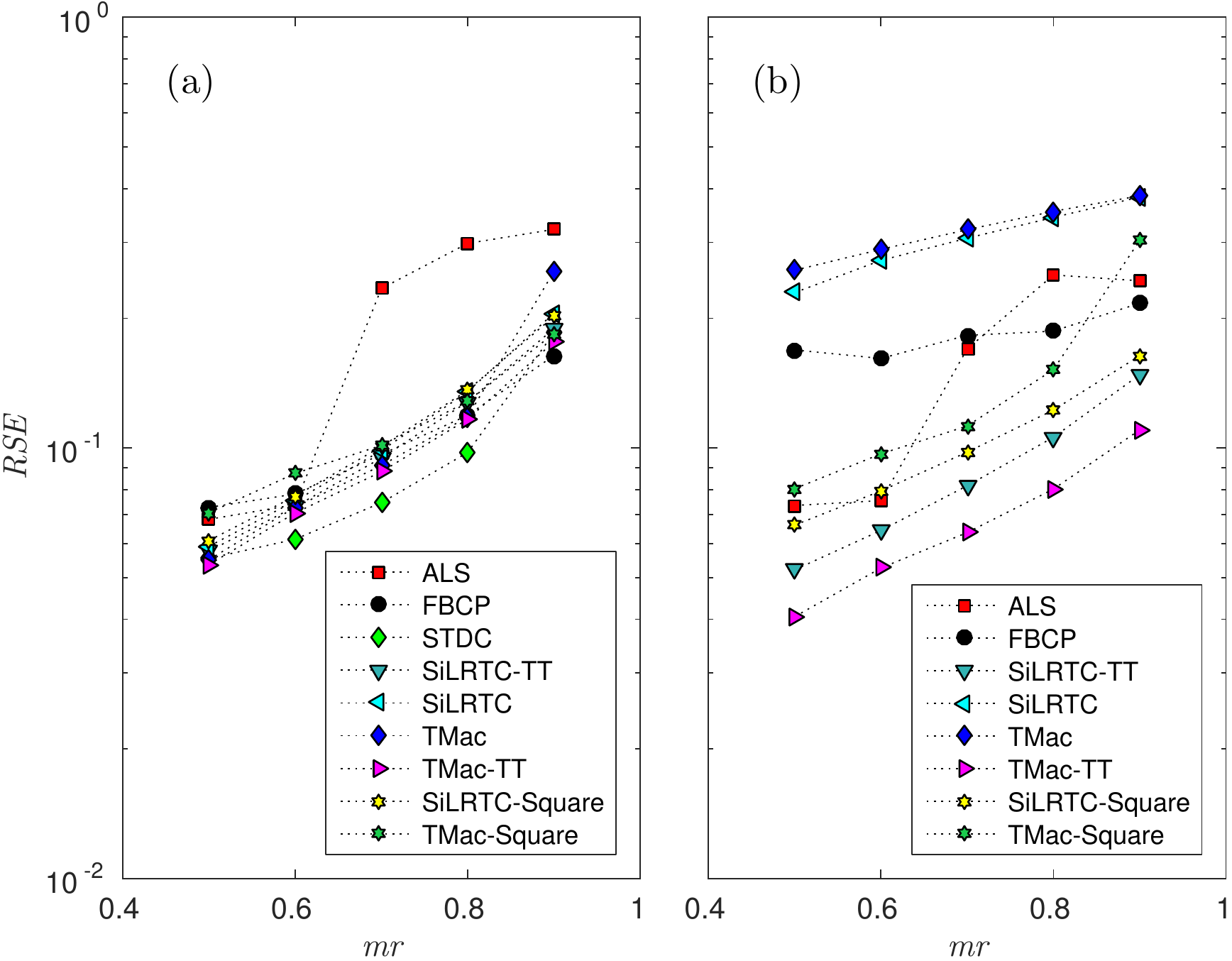}\\
	\caption{Performance comparison between different tensor completion algorithms based on the RSE vs the missing rate when applied to the Lena image. (a) Original tensor (no order augmentation). (b) Augmented tensor using KA scheme.}
	\label{fig7}
\end{figure}

The results for the experiment performed on the Lena image and recovery results for $mr=0.9$ are shown in Fig.~\ref{fig8} and Fig.~\ref{fig7}, respectively. The results show that TMac-TT gives the best results (lowest $RSE$) for each $mr$ when using the KA scheme.

For the House image, the missing entries are now chosen as white text, and hence the missing rate is fixed. The result is shown in Fig.~\ref{fig9}. STDC provides the best performance without augmentation, while all other algorithms are comparable. However, the outlines of text can still be clearly seen on the STDC image. Using tensor augmentation, TMac-TT and TMac-Square provides the best performance, where the text is almost completely removed using TMac-TT.


\subsection{Video completion with ket augmentation}
In color video completion we benchmark FBCP, ALS, TMac, TMac-TT and TMac-Square against two videos, \textit{New York City (NYC)} and \textit{bus}\footnote{Videos available at {https://engineering.purdue.edu/\textasciitilde reibman/ece634/}}.  The other methods are computationally intractable or not applicable for $N\geq4$ in this experiment. For each video, the following preprocessing is performed: Resize the video to a tensor of size $81\times729\times1024\times3$ ($frame\times image\ row\times image\ column\times RGB$). The first frame of each video can be seen in Figs. \ref{fig:a} and \ref{fig:b}.  The $frame$ mode is merged with the $image\ row$ mode to form a third-order tensor, which we define here as a \textit{video sequence tensor (VST)}, of size $59,049\times1024\times3$ ($combined\ row\times image\ column\times RGB$). Examples of the VST can be seen in the range 20000:20700 for $combined\ row$ in Figs. \ref{fig:c} and \ref{fig:d}. Hence, rather than performing an image completion on each frame, we perform our tensor completion benchmark on the \textit{entire video}. It is important to highlight that we only benchmark with a ket augmented (not the third-order) VST due to computational intractability for high-dimensional low-order tensors.
		
		\begin{figure*}[htpb!]
			\centering
			\includegraphics[scale=0.55]{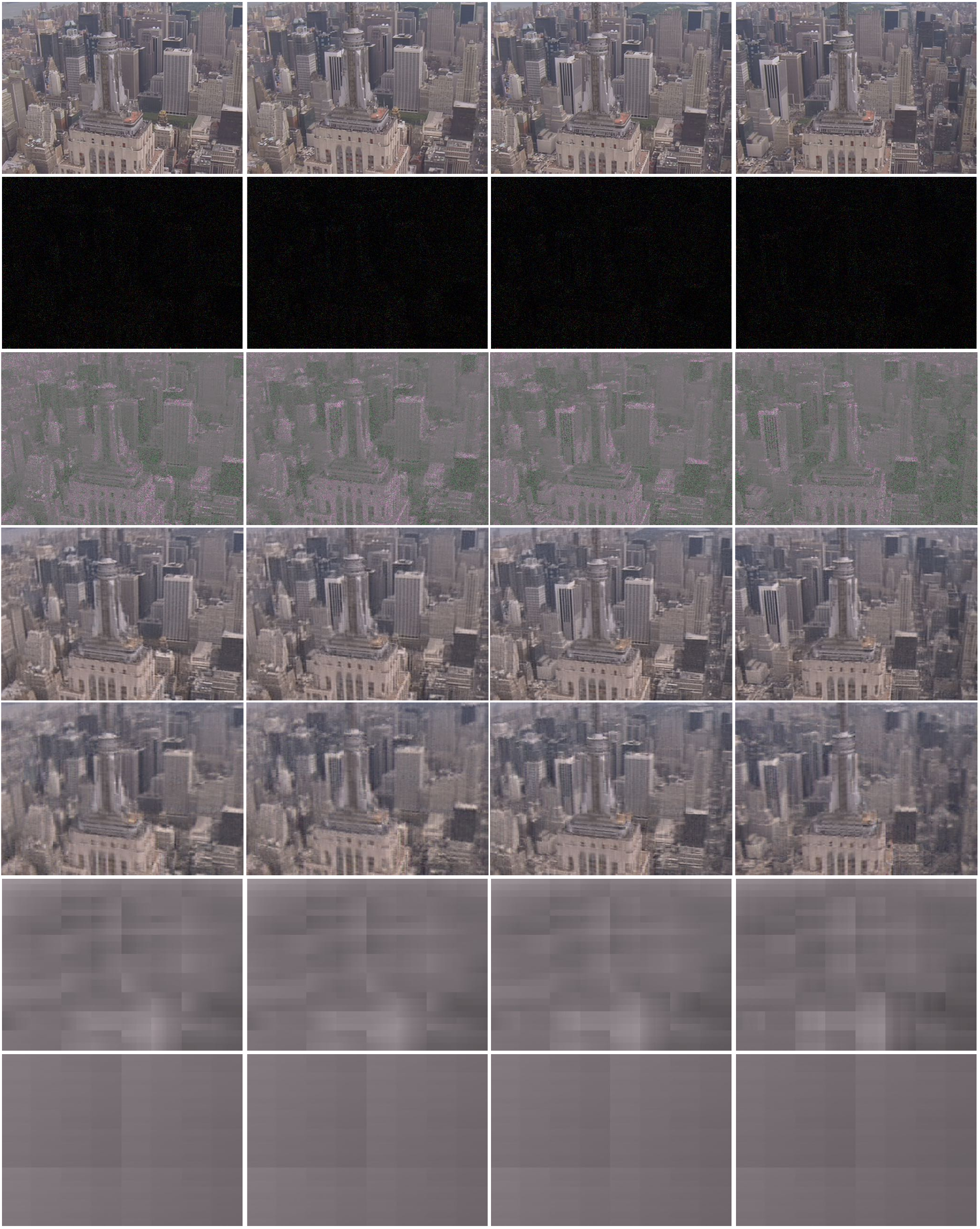}
			\caption{The 7th, 21st, 33rd and 70th frames (from left to right column) in the NYC video, with each row (from top to bottom) representing the original frames, original frames with 95\% missing entries, TMac, TMac-TT, TMac-Square, ALS and FBCP.}\label{citycompare}
		\end{figure*}	
			
		\begin{figure*}[htpb!]
			\centering
			\includegraphics[scale=0.55]{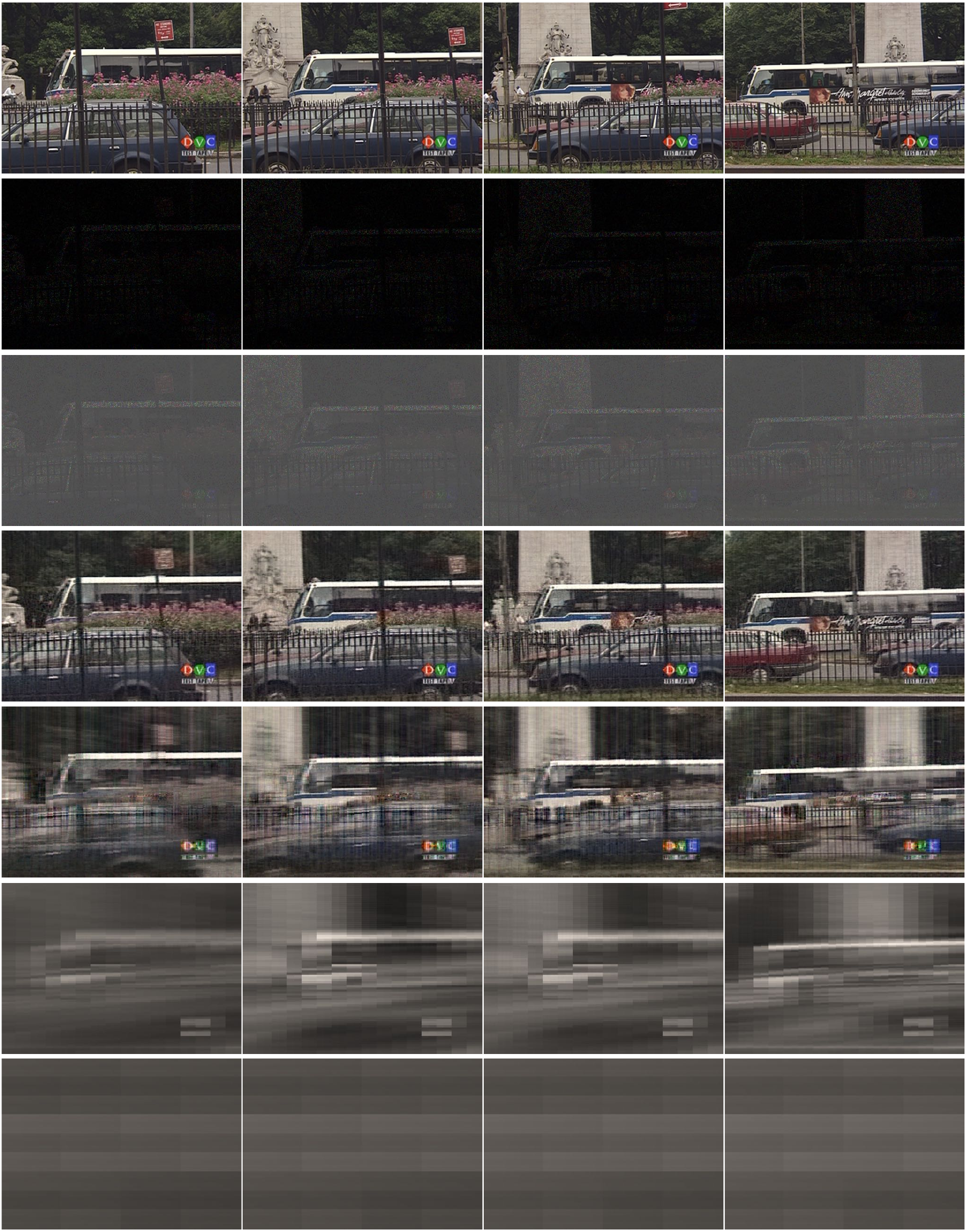}
			\caption{The 7th, 21st, 33rd and 70th frames (from left to right column) in the bus video, with each row (from top to bottom) representing the original frames, original frames with 95\% missing entries, TMac, TMac-TT, TMac-Square, ALS and FBCP.}\label{buscompare}
		\end{figure*}
			
Using KA, reshape the VST to a low-dimensional high-order tensor of size $6\times6\times6\times6\times6\times6\times6\times6\times6\times6\times3$. The eleventh-order VST is directly used for the tensor completion algorithms.
\begin{table}[!htbp]
	\caption{RSE and SSIM tensor completion results for 95\%, 90\% and 70\% missing entries from the NYC video.}
	\label{tableNYC}
	\centering 
	\begin{tabular}{c|cc|cc|cc}
		  & \multicolumn{2}{c|}{$mr=0.95$} & \multicolumn{2}{c|}{$mr=0.9$} & \multicolumn{2}{c}{$mr=0.7$}  \\
		\hline
		Algorithm & RSE & SSIM & RSE &  SSIM & RSE & SSIM\\				
		\hline
		{FBCP}&0.210&0.395&0.210&0395&0.210&0.396\\
		{ALS}&0.193&0.397&0.189&0.398&0.168&0.429\\
		{TMac}&0.185&0.605&0.143&0.750&0.055&\bf{0.967}\\
		{TMac-TT}&\bf{0.072}&\bf{0.876}&{\bf 0.066}&{\bf 0.902}&\bf{0.053}&0.949\\	
		{TMac-Square}&0.111&0.722&0.076&0.901&0.056&0.946\\		
		\hline		
	\end{tabular}
\end{table}

For the case of 95\% missing entries, results of the benchmark can be seen in Figs. \ref{citycompare} and \ref{buscompare}. The NYC results in Fig. \ref{citycompare} shows that FBCP and ALS are completely incomprehensible, whereas only the TMac-based algorithms can successfully complete the video. Moreover, in this case, TMac-TT outperforms all algorithms, which can be seen with the $RSE$ and mean structural similarity index ($SSIM$) \cite{1284395} (over all 81 \textit{frames}) in Table \ref{tableNYC} for $mr=0.95$. For the bus results in Fig. \ref{buscompare}, TMac-TT outperforms all algorithms. The other TT rank-based algorithm ALS can only manage a simple structure of the bus, and FBCP cannot produce any resemblence to the original video. Table \ref{tableBUS} summarizes the $RSE$ and mean $SSIM$ results. With 90\% missing entries, the results are similar to those of 95\% missing entries, however, TMac-TT and TMac-Square are now comparable in performance for the NYC video.
For the NYC video with $mr=0.7$, Table \ref{tableNYC} shows that all TMac-based algorithms are comparable, with FBCP and ALS unable to reproduce a sufficient approximation. In the bus video, TMac-TT and TMac-Square provide comparable RSE and SSIM.
\begin{table}[!htbp]
	\caption{RSE and SSIM tensor completion results for 95\%, 90\% and 70\% missing entries from the bus video.}
	\label{tableBUS}
	\centering 
	\begin{tabular}{c|cc|cc|cc}
		 & \multicolumn{2}{c|}{$mr=0.95$} & \multicolumn{2}{c|}{$mr=0.9$} & \multicolumn{2}{c}{$mr=0.7$}  \\
		\hline
		Algorithm & RSE & SSIM & RSE &  SSIM & RSE & SSIM\\				
		\hline
		{FBCP}&0.527&0.269&0.527&0.269&0.504&0.271\\
		{ALS}&0.447&0.323&0.342&0.387&0.271&0.513\\
		{TMac}&0.518&0.316&0.496&0.374&0.402&0.598\\
		{\bf TMac-TT}&\bf{0.154}&\bf{0.807}&{\bf 0.092}&{\bf 0.932}&\bf{0.062}&\bf{0.974}\\	
		{TMac-Square}&0.267&0.582&0.196&0.781&0.077&0.968\\		
		\hline		
	\end{tabular}
\end{table}	

In summary, the bus video includes more vibrant colours and textures compared to the NYC video, which can be clearly seen from the overall $SSIM$ performance in Tables \ref{tableNYC} and \ref{tableBUS}. It is important to highlight that TMac-TT still provides a high quality ($SSIM=0.807$) approximation for the high missing ratio ($mr=0.95$) test, where the next best result of TMac-Square had only $SSIM=0.582$. This demonstrates the superiority of using TMac-TT over the other algorithms for high missing ratio video completion problems.			
\section{Conclusion \label{sec5}}
{\color{black}A novel approach} to the LRTC problem {\color{black}based on} TT rank was introduced along with
corresponding algorithms for its solution. The SiLRTC-TT algorithm was defined to minimize the TT rank of a tensor by TT nuclear norm optimization. Meanwhile, TMac-TT was proposed, which is based on the multilinear matrix factorization model to minimize the TT-rank. The latter is more computationally efficient due to the fact that it does not need the SVD. The proposed algorithms are employed to simulate both synthetic and real world data represented by higher-order tensors. For synthetic data,  our algorithms prevails against the others when the tensors have low TT rank. Their performance is comparable in the case of low Tucker rank tensors. The TT-based algorithms are quite promising and reliable when applied to real world data. To validate this, we studied image and video completion problems.
Benchmark results show that when applied to original tensors without tensor augmentation, our algorithms are comparable to STDC in image completion. However, in the case of augmented tensors, our proposed algorithms not only outperform the others, but also provide better recovery results when compared to the case without tensor order augmentation in both image and video completion experiments.

{\color{black} Applications of the proposed TT rank optimization based tensor completion
to data compression, text mining, image classification and video indexing are under our interest.}

\bibliographystyle{IEEEtran}
\bibliography{IEEEabrv,TC}
\end{document}